\documentclass[11pt,preqno,a4paper]{amsart}

\usepackage{enumitem}
\usepackage{amssymb,amsmath,mathtools}
\usepackage{hyperref}

\usepackage{microtype}
\usepackage{bbm} 
\usepackage{a4wide}
\newcommand*{\mailto}[1]{\href{mailto:#1}{\nolinkurl{#1}}}

 \allowdisplaybreaks
 
\newtheorem{theorem}{Theorem}[section]
\newtheorem{definition}{Definition}[section]
\newtheorem{lemma}[theorem]{Lemma}
\newtheorem{proposition}[theorem]{Proposition}
\newtheorem{corollary}[theorem]{Corollary}
\newtheorem{remark}[theorem]{Remark}

\newcommand{\be}{\begin{equation}}
\newcommand{\ee}{\end{equation}}

\numberwithin{equation}{section}

\makeatletter
\newcommand{\dlmf}[1]{%
	\cite[%
	\def\nextitem{\def\nextitem{, }}%
	\@for \el:=#1\do{\nextitem\href{http://dlmf.nist.gov/\el}{(\el)}}%
	]{dlmf}%
}
\makeatother

\begin{document}

\allowdisplaybreaks

\title [Asymptotic stability of Landau solutions]{Asymptotic stability of Landau solutions to Navier-Stokes system under $L^p$-perturbations}
\author{Yanyan Li         \and
        Jingjing Zhang \and
        Ting Zhang*
}

\address{Yanyan Li,  Department of Mathematics, Rutgers University \\ NJ 08854, USA
}
\email{yyli@math.rutgers.edu }
\address{Jingjing Zhang and Ting Zhang*, School of Mathematical Sciences, Zhejiang University\\
 Hangzhou 310027, China
}
\email{zhangjingjing94@126.com;  zhangting79@zju.edu.cn  }


\begin{abstract}
In this paper, we show that Landau solutions to the Navier-Stokes system are asymptotically stable under $L^3$-perturbations. We give the local well-posedness of solutions to the perturbed system with  initial data in $L_{\sigma}^3$ space and the global  well-posedness with  small initial data in $L_{\sigma}^3$ space, together with a study of the $L^q$ decay for all $q>3.$ Moreover, we have also studied the local well-posedness, global  well-posedness and stability in $L^p$ spaces for $3<p<\infty$.
\end{abstract}

\maketitle

\section{Introduction}
\label{intro}

The Cauchy problem for the incompressible Navier-Stokes system in $\{(x,t)| x\in \mathbb{R}^3,t\geq 0 \} $ with given initial data and external force has the form
\begin{equation}\label{NS}
\begin{cases}
 u_{t}-\Delta u+(u \cdot \nabla) u+\nabla p =f, \\
 \nabla\cdot u =0, \\
 u(x, 0) =u_{0}(x),
\end{cases}
\end{equation}
where $u=(u_1,u_2,u_3)$ and  $p$ denote the velocity field and pressure respectively.

Note that when we consider the construction of solutions to the Cauchy problem (\ref{NS}), there are essentially two methods: the energy method and the perturbation theory. The energy method is based on $a$-$priori$ energy estimate
\begin{eqnarray*}
\int_{\mathbb{R}^3}|u(x,t)|^2dx+\int_0^t\int_{\mathbb{R}^3}2|\nabla u(x,s)|^2dxds
&\leq& \int_{\mathbb{R}^3}|u_0(x)|^2dx+\int_0^t\int_{\mathbb{R}^3}2(f,u)(x,s)dxds.
\end{eqnarray*}
The global existence of weak solutions was established by Leray \cite{Ler} for divergence free initial data  $u_0\in L^2(\mathbb{R}^3)$ and $f=0$. The energy method gives the existence, but the uniqueness and regularity for  solutions still remain open, see $e.g.$ \cite{Ba,Caf,Can1,Jia,Lem1,Lem,Te} and references therein.

As for the perturbation theory, we treat the nonlinear term $(u \cdot \nabla) u$ as a perturbation and use the scaling property to choose function spaces. As we know,  system (\ref{NS}) has the natural scaling
\begin{equation*}
\begin{aligned}
u_{\lambda}(x,t)=\lambda u(\lambda x,\lambda^2 t),\quad
p_{\lambda}(x,t)=\lambda ^2p(\lambda x,\lambda^2 t).
\end{aligned}
\end{equation*}
Therefore, the space $L^3(\mathbb{R}^3)$ is a well-known simple example of scaling-invariant space. By the Duhamel principle, we can write these solutions into an integral formulation
\begin{equation}
u(x,t)=e^{t\Delta}u_0+\int_0^t e^{(t-s)\Delta}\mathbb{P}(f-u\cdot\nabla u)ds,
\end{equation}
where $\mathbb{P}$ denotes the Leray projector which projects on divergence-free vector fields.
Solutions constructed in this way are called mild solutions. Usually, by means of the contraction mapping principle, we can obtain global well-posedness of mild solutions to system (\ref{NS})  with small enough initial data in appropriate scaling-invariant spaces. We refer readers to \cite{Can0,Can1,Kat,Ko,Lem1,Lem,Rob}  for additional background and references.

There are many results on  the existence of weak solutions and $L^2$-decay of weak solutions of the Navier-Stokes system, see e.g. \cite{Kat,Ma,Sch,Wie}  and references therein.  When $f=0,$ the $L^2$-decay of weak solutions to system (\ref{NS}) can be viewed as the global asymptotic stability in $L^2$ of the trivial solution $(u,p)=(0,0).$ Later, Borchers and Miyakawa \cite{Bo} addressed similar questions on the global asymptotic stability to a family of stationary solutions.

The stationary Navier-Stokes system in $ \mathbb{R}^3 $ has the form
\begin{equation}\label{SNS}
\begin{cases}
 -\Delta v+(v \cdot \nabla) v+\nabla p =f, \\
 \nabla\cdot v =0. \\
\end{cases}
\end{equation}
When $f=(b(c)\delta_0,0,0)$ with $b(c)=\frac{8\pi c}{3(c^2-1)}\left(2+6c^2-3c(c^2-1)\ln\left(\frac{c+1}{c-1}\right)\right)$ and $\delta_0$ the Dirac measure, $(v_c, p_c)$ given by the following formulas
\begin{equation}\label{v_c p_c}
\begin{aligned}{v_{c}^{1}(x)=2 \frac{c|x|^{2}-2 x_{1}|x|+c x_{1}^{2}}{|x|\left(c|x|-x_{1}\right)^{2}},} & \ {v_{c}^{2}(x)=2 \frac{x_{2}\left(c x_{1}-|x|\right)}{|x|\left(c|x|-x_{1}\right)^{2}}}, \\ {v_{c}^{3}(x)=2 \frac{x_{3}\left(c x_{1}-|x|\right)}{|x|\left(c|x|-x_{1}\right)^{2}},} & \ {p_{c}(x)=4 \frac{c x_{1}-|x|}{|x|\left(c|x|-x_{1}\right)^{2}}},\end{aligned}
\end{equation}
with $|x|=\sqrt{x_{1}^{2}+x_{2}^{2}+x_{3}^{2}}$ and constant  $|c|>1$ are distributional solutions to system (\ref{SNS}) in $\mathbb{R}^3$. We note that $b(c)$ is decreasing on $(-\infty,-1)$ and $(1, \infty)$, ${\lim}_{c\rightarrow 1}b(c)=+\infty,$ ${\lim}_{c\rightarrow -1}b(c)=-\infty$ and ${\lim}_{|c|\rightarrow \infty}b(c)=0.$ The explicit stationary solutions (\ref{v_c p_c}) were discovered   by  Landau \cite{La}. These solutions have been called  Landau solutions. Tian and Xin \cite{Ti} proved  that all  $(-1)-$homogeneous, axisymmetric nonzero
solutions of system (\ref{SNS}) in $C^2(\mathbb{R}^3\backslash \{0\})$ are Landau solutions.
 $\check {\mathrm S}$ver$\acute {\mathrm a}$k \cite{Sv} proved that Landau solutions are the only $(-1)-$homogeneous solutions  in $C^2(\mathbb{R}^3\backslash \{0\})$. More details can be found in \cite{Can,La,Sl,Sv,Ti}.

Karch and Pilarczyk \cite{Gr} showed that Landau solutions are asymptotically stable under any $L^2$-perturbations. The crucial role played in their paper is an application of the
Hardy-type inequality
\begin{equation}
\left|\int_{\mathbb{R}^3}w\cdot (w \cdot \nabla)v_c dx\right|\leqq K(c)\|\nabla \otimes w\|_2^2,
\end{equation}
where positive function $K(c)=12\max _{j, k \in\{1,2,3\}} K_{j, k}(c)$ with  functions $K_{j, k}:(-\infty,-1) \cup(1, \infty) \rightarrow$
$(0, \infty)$ for every $j, k \in\{1,2,3\}$ satisfying
\begin{equation}\label{K_{j,i}}
\left|\partial_{x_{j}} v_{c}^{k}(x)\right| \leqq \frac{K_{j, k}(c)}{|x|^{2}}.
\end{equation}
Moreover, $K_{j,k}(c)$ satisfies
\begin{equation} \label{about c}
\lim_{|c|\rightarrow 1}K_{j,k}(c)=+\infty  \text{  and  } \lim_{|c|\rightarrow \infty}K_{j,k}(c)=0,
\end{equation}

In 2017, Karch, Pilarczyk and Schonbek \cite{Ka} generalized the work of \cite{Gr}. They gave a new method to show the $L^2$-asymptotic stability of a large class of global-in-time solutions including the Landau solutions. Their work also generalizes results in a series of articles on $L^2$-asymptotic stability either of the zero solution \cite{Bor1,Kaj,Oga,Sch1,Sch,Wie} or nontrivial stationary solutions \cite{Bo} to system (\ref{NS}).
The above results give the existence  in $L^2$ space.  The uniqueness in $L^2$ space is a major open problem. We will consider the stability of Landau solutions to the Navier-Stokes system  in $L^p$ spaces with $3\leq p <\infty$.

We denote  $(u,p)(x,t)$  solution to the Navier-Stokes system (\ref{NS}) with the given external force $f=(b(c)\delta_0,0,0)$ and initial data $u_0=v_c+ w_0$. By a direct calculation, functions $w(x,t)=u(x,t)-v_c(x)$ and $\pi(x,t)=p(x,t)-p_c(x)$ satisfy the following system
\begin{equation}\label{PNS}
\begin{cases}
 w_{t}-\Delta w+(w \cdot \nabla) w+(w \cdot \nabla) v_{c}+\left(v_{c} \cdot \nabla\right) w+\nabla \pi=0,\\
 \nabla\cdot w=0,\\
 w(x, 0) =w_{0}(x).
\end{cases}
\end{equation}
  We will consider the well-posedness of solutions to system (\ref{PNS}) in $L^p$ spaces with $3\leq p <\infty$. We can obtain  global  well-posedness  of solutions to  system  (\ref{PNS}) with  small initial data in $L_{\sigma}^3$ space and local well-posedness  with general initial data in $L^3_{\sigma}$ space, see Theorem \ref{L^3 well-posedness}. For initial data $w_0\in L^p_{\sigma}$ with $3<p<\infty$, we have  local well-posedness results, see Theorem \ref{p>3 result}. In addition, for general initial data in $L^3_{\sigma},$ we have the global existence of $L^2+L^3$ weak solutions, see  Definition \ref{global 3} and Theorem \ref{w-global-weak-existence}.

Karch and Pilarczyk \cite{Gr} denote the linear operator $\mathcal{L}$
\begin{equation}\label{L def}
 \mathcal{L}u=-\Delta u+\mathbb{P}\left((u \cdot \nabla) v_{c}\right)+\mathbb{P}\left(\left(v_{c} \cdot \nabla\right) u\right),
\end{equation}
where $\mathbb{P}$ is the Leray projector.
For system (\ref{PNS}), we can write  solution in the following formula
\begin{equation}\label{eq}
w(x,t)=e^{-t\mathcal{L}}w_0-\int_0^te^{-(t-s)\mathcal{L}}\mathbb{P}\nabla\cdot(w\otimes w)ds:= a+N(w,w).
\end{equation}

Karch and Pilarczyk \cite{Gr} showed that  $-\mathcal{L}$ is the infinitesimal generator of an analytic  semigroup of bounded linear operators on $L_{\sigma}^2({\mathbb{R}^{3}})$. We show that for $1< q<\infty,$  $-\mathcal{L}$ is the infinitesimal generator of an analytic  semigroup of bounded linear operators on $L_{\sigma}^q({\mathbb{R}^{3}})$, see Theorem \ref{linear ope property} in Section \ref{linear}.

\subsection{$L^p$  mild solutions, $3\leq p<\infty$}

Let us give the following standard definition of $L^p$  mild solutions, $3\leq p<\infty$.
\begin{definition}
Let $3\leq p<\infty$ and $T>0,$ a function $w$ is a $L^p$ mild solution of system (\ref{PNS}) with initial data $w_0\in L_\sigma^p(\mathbb{R}^3)$  on $[0,T]$ if
\begin{equation}\label{w_space}
w\in C([0,T];L_\sigma^p(\mathbb{R}^3))\cap L^{\frac{4p}{3}}([0,T];L_\sigma^{2p}(\mathbb{R}^3)),
\end{equation}
and
\begin{equation}\label{w_formula}
w(x,t)=e^{-t\mathcal{L}}w_0-\int_0^te^{-(t-s)\mathcal{L}}\mathbb{P}\nabla\cdot(w\otimes w)ds.
\end{equation}
\end{definition}

This  solution is  global if (\ref{w_space}) and (\ref{w_formula}) hold for any $0<T<\infty.$\

In the above, $e^{-t\mathcal{L}}$ denotes the analytic semigroup of bounded linear operators on  $L_{\sigma}^p(\mathbb{R}^3)$ generated by $-\mathcal{L}$. See Lemma \ref{alemm}, Lemma \ref{p_alemm} and Theorem \ref{linear ope property}. Properties of $\int_0^te^{-(t-s)\mathcal{L}}\mathbb{P}\nabla\cdot(w\otimes w)ds$ can be seen in Lemma \ref{lem2}.

Now we give the following theorem which shows the well-posedness results in $L^3$ space and the decay rates of solutions to system (\ref{PNS}).
\begin{theorem}\label{L^3 well-posedness}
There exist positive universal constants $c_3,$ $\varepsilon_0$ and $C$ with the following properties\\
$(i)$ For every $|c|>c_3$ and $w_0\in L_{\sigma}^3(\mathbb{R}^3)$, there exists a positive constant $T$ depending only on $w_0$ such that system (\ref{PNS}) has a unique $L^3$ mild solution $w$ on $[0,T]$.  Moreover,   $\nabla (|w|^{\frac32})\in L^2([0,T];L^2(\mathbb{R}^3)).$\\
$(ii)$ If in addition, $\|w_0\|_{L^3(\mathbb{R}^3)}<\varepsilon_0$, then  system (\ref{PNS}) has a  unique global $L^3$ mild solution $w$. Moreover,  $\nabla (|w|^{\frac32})\in L^2([0,\infty);L^2(\mathbb{R}^3))$,
\begin{equation}\label{L^3 stability-1}
\|w\|_{C_t(L_x^3)\cap L^4_t(L^6_x)} + \|\nabla(|w|^{\frac32})\|^{\frac23}_{L^2_tL^2_x}\leq C \|w_0\|_{L^3(\mathbb{R}^3)},
\end{equation}
and
\begin{equation}\label{L^3 stability-2}
\lim_{t\rightarrow \infty}\|w(t)\|_{L^3(\mathbb{R}^3)}=0.
\end{equation}
$(iii)$ For any $q>3$, there exists a  positive constant $\tilde{c}_q$ depending only on $q$ such that when
$|c|>\tilde{c}_q$, the solution in $(ii)$ satisfies
$$
w\in L^{\infty}([\tau,\infty), L_\sigma^q(\mathbb{R}^3)), \quad \text{for all} \quad   \tau >0,
$$
and
\begin{equation}\label{nonlinear decay}
\|w(t)\|_{L^q(\mathbb{R}^3)} \leq (\frac13-\frac1q)^{\frac32 (\frac13-\frac1q)} t^{\frac{3}{2q}-\frac{1}{2}}\|w_0\|_{L^3(\mathbb{R}^3)},\quad \text{for all} \quad  t>0.
\end{equation}

\end{theorem}

\begin{remark}
From  (\ref{|c|2}), (\ref{|c|3}) and (\ref{mu}) in this paper, we see a more detailed dependence of $\tilde{c}_q$. On the other hand, we tend to believe that $c$ can be chosen as a constant independent of $q$, and we plan to investigate $L^\infty$ decay in our future work.
\end{remark}

\begin{remark}
 It follows from Theorem \ref{L^3 well-posedness}  that the flow described by the Landau solution is asymptotically stable under   $L^3$-perturbations.
\end{remark}

\begin{remark}
For the two-dimensional Navier-Stokes system,  Carlen and Loss \cite{Ca} gave the decay rate of solutions to the vorticity equation. We adapt the method in  \cite{Ca} to give the decay rate of solutions to  system (\ref{PNS}), and we treat  the pressure term $\pi$ by using the $A_p$ weight inequalities  for the Riesz transforms \cite{Gra,St}.
\end{remark}
For $3<p<\infty,$ we have the following result
\begin{theorem}\label{p>3 result}
For $p\in (3,\infty)$ and  $w_0\in L_{\sigma}^p(\mathbb{R}^3)$, there exist two constants ${c}_p$ and $T$, where ${c}_p$ depends only on $p$ while $T$ depends only on $p$ and an upper bound of $\|w_0\|_{L^p}$ such that for all  $|c|>{c}_p,$  system (\ref{PNS}) has  a unique $L^p$ mild solution $w$ on $[0,T]$. Moreover,   $\nabla(|w|^{\frac p2}) \in L_t^2([0,T];L_x^2(\mathbb{R}^3)).$ If in addition, $w_0\in L_{\sigma}^p\cap L_{\sigma}^3(\mathbb{R}^3)$, $\|w_0\|_{L^3}<\varepsilon_0$, where $\varepsilon_0$ is as in Theorem \ref{L^3 well-posedness}, there exists a  unique global $L^p$ mild solution $w$ to system (\ref{PNS}). Moreover, for some universal constant $C,$
\begin{equation}\label{p>3-1}
\|w \|_{C_t(L_x^p)\cap L_t^\frac{4p}{3}(L_x^{2p})}+\|\nabla(|w|^{\frac p2})\|_{L_t^2L_x^2}^{\frac2p} \leq C \|w_0\|_{L^p}.
\end{equation}
\end{theorem}
\begin{remark}
Under the condition of Theorem \ref{p>3 result}, applying  Theorem \ref{L^3 well-posedness},  (\ref{L^3 stability-1})-(\ref{L^3 stability-2}) hold. And, if in addition, $|c|>\tilde{c}_q,$ then according to Theorem \ref{L^3 well-posedness}, (\ref{nonlinear decay}) holds.
\end{remark}

\begin{remark}
Note that $w_0\in L^p $ with $p>3$ implies $w_0\in L^2_{uloc}$. For the Navier-Stokes system with $u_0\in L^2_{uloc}$, several authors \cite{Bas,Kik,Lem1,Lem} gave the local existence of weak solution $u$. Moreover, global weak solution exists for decaying initial data $u_0\in E_2$ with
  $$E_2=\left\{f \in L_{\text {uloc }}^{2}:\|f\|_{L^{2}\left(B\left(x_{0}, 1\right)\right)} \rightarrow 0,  \text { as }\left|x_{0}\right| \rightarrow \infty\right\}.$$  
  Kown and Tsai \cite{Kw} generalized the  global existence with non-decaying initial data whose local oscillations decay. Very recently, J.J.  Zhang and T. Zhang \cite{Zhang}  have given the local existence of solutions to system (\ref{PNS}) with initial data $w_0\in L^p_{uloc}$, $p\geq 2.$
In view of  these results, we plan  to study the global existence of weak solutions to system (\ref{PNS}) with  initial  data $w_0\in L_{\sigma}^p$ for $p>3$ in our future work.
\end{remark}
\begin{remark}
L. Li, Y.Y. Li and X. Yan investigated homogeneous solutions of stationary Navier-Stokes
system with isolated singularities on the unit sphere \cite{Li1,Li2,Li4,Li3}.  For a subclass of (-1)-homogeneous axisymmetric no-swirl solutions on the unit sphere minus north and south poles classified in \cite{Li2},  Y.Y. Li and X. Yan have proved in \cite{Yan Li} the asymptotic stability under $L^2$-perturbations.  We will focus on these  homogeneous solutions  in our future work.
\end{remark}

Results in Theorems \ref{L^3 well-posedness} and \ref{p>3 result} show the existence and uniqueness of solution $w$ to system (\ref{PNS}) in corresponding space. Actually the solution depends continuously on initial data.

\begin{theorem}\label{continuty}
For every $|c|>c_p$ and  $u_0\in L_{\sigma}^p(\mathbb{R}^3)$ with $3\leq p<\infty,$ assume that $u$ is the unique mild solution to system (\ref{PNS}) on $[0, T_{max})$. Then, for any $T\in (0, T_{max}),$ there exists $\varepsilon>0$ such that for any  $v_0\in L_{\sigma}^p(\mathbb{R}^3)$, $\|u_0-v_0\|_{L^p}<\varepsilon,$  there exists a unique  $L^p$ mild solution $v$ on $[0,T]$ with initial data $v\vert_{t=0}=v_0$. Moreover,
\begin{equation}\label{sol}
\lim_{u_0\rightarrow v_0 \text{\ in\ }L^p}\left(\|u-v\|_{C_TL_x^p\cap L^{\frac{4p}{3}}_TL_x^{2p} }+\left\|\nabla\left(|u-v|^{\frac{p}{2}}\right)\right\|^{\frac2p}_{L_T^2L_x^2}\right)= 0.
\end{equation}
\end{theorem}

The constant $c_p$ in the above theorem is the one given in Theorems \ref{L^3 well-posedness} and \ref{p>3 result}.

\begin{remark}
Karch, Pilarczyk and Schonbek \cite{Ka}  showed the $L^2$-asymptotic stability of a large class of global-in-time solutions including the Landau solutions. Based on similar  proof of Theorem \ref{continuty}, we can obtain the $L^3$-asymptotic stability of a  class of solutions $v_c+w$, where $w$ is as in Theorem \ref{L^3 well-posedness}.  More precisely,  letting $V$ as perturbation of $v_c+w$, when $\|V_0\|_{L^3}\leq (4C^2e^{2C \int_0^{\infty} \|w\|_{L^6}^4 dt})^{-1}$, using the similar  proof of  (\ref{Z-crucial}), we obtain
\begin{eqnarray}
\|V\|_{C_{t}L_x^3\cap L^4_{t}L_x^6 }+\left\|\nabla\left(|V|^{\frac{3}{2}}\right)\right\|^{\frac23}_{L_{t}^2L_x^2}\leq  2C\|V_0\|_{L^3}e^{C\int_0^{\infty} \|w\|_{L^6}^4 dt}.
\end{eqnarray}

 \end{remark}

 \subsection{Weak solution}

 Karch and Pilarczyk \cite{Gr,Ka} proved the following results: for every $w_{0} \in L_{\sigma}^{2}\left(\mathbb{R}^{3}\right)$, there exists a global weak solution
 $$w\in C_{w}\left([0, T], L_{\sigma}^{2}\left(\mathbb{R}^{3}\right)\right) \cap L^{2}\left([0, T], \dot{H}_{\sigma}^{1}\left(\mathbb{R}^{3}\right)\right)$$
  for every $T>0$ which satisfies the strong energy inequality
\begin{equation}\label{energy-inequality}
\|w(t)\|_{2}^{2}+2(1-K(c)) \int_{s}^{t}\|\nabla \otimes w(\tau)\|_{2}^{2} \mathrm{d} \tau \leqq\|w(s)\|_{2}^{2}
\end{equation}
for almost all $s \geq 0,$ including $s=0$ and all $t \geq s$. The definition of the weak solution is as follows
\begin{definition} ($L^2$-weak solution)
For $w_0\in L_{\sigma}^2(\mathbb{R}^3)$, a function $w$ is a $L^2$-weak solution of system (\ref{PNS}) on $[0,T]$ if\\
i) $w\in C_w([0,T];L_{\sigma}^2(\mathbb{R}^3))\cap L^2 ([0,T];\dot{H}_{\sigma}^1(\mathbb{R}^3)).$\\
ii)
For all $t \geqq s \geqq 0$, all $\varphi \in C\left([0, \infty), H_{\sigma}^{1}\left(\mathbb{R}^{3}\right)\right) \cap C^{1}\left([0, \infty), L_{\sigma}^{2}\left(\mathbb{R}^{3}\right)\right),$
\begin{eqnarray*}
&&(w(t), \varphi(t))+\int_{s}^{t}\left[(\nabla w, \nabla \varphi)+(w \cdot \nabla w, \varphi)+\left(w \cdot \nabla v_{c}, \varphi\right)+\left(v_{c} \cdot \nabla w, \varphi\right)\right] \mathrm{d} \tau\\
&=&(w(s), \varphi(s))+\int_{s}^{t}\left(w, \varphi_{\tau}\right) \mathrm{d} \tau.
\end{eqnarray*}
iii) For all $\phi\in C_{c}^{\infty}(\mathbb{R}^{3})$, $\lim_{t\rightarrow0}\int_{\mathbb{R}^3}w\cdot \phi dx=\int_{\mathbb{R}^3}w_0\cdot \phi(0) dx$.\\
iv) $w$ satisfy the energy inequality
\begin{eqnarray*}
&&\int_{\mathbb{R}^{3}}|w|^{2} \xi(x, t) d x+2 \int_{0}^{t} \int_{\mathbb{R}^{3}}|\nabla w|^{2} \xi dx ds\\
 &\leq &\int_{\mathbb{R}^{3}}|w_0|^{2} \xi(x, 0) d x+\int_{0}^{t} \int_{\mathbb{R}^{3}}(2v_c\otimes w:\nabla w \xi \\
&&+ \left(\partial_{s} \xi+\Delta \xi\right)|w|^{2}+(|w|^{2}+2 \pi+2v_c\cdot w)(w \cdot \nabla) \xi + |w|^2v_c\cdot \nabla \xi )dx ds,
\end{eqnarray*}
for any $t \in[0, T]$  and for all non-negative smooth functions $\xi \in C_c^{\infty}([0,T]\times \mathbb{R}^{3}).$\\
\end{definition}

The following is a weak-strong uniqueness theorem which is analogous to the one for the Navier-Stokes system (Theorem 4.4 in \cite{Tsa}).
\begin{theorem}\label{s-w-uniqueness}
Let $|c|>8\sqrt2 +1$, $w_0\in L_{\sigma}^2(\mathbb{R}^3)$. Assume that $u, v$ are $L^2$-weak solutions of system (\ref{PNS}) on $[0,T]$ with initial data $u\vert_{t=0}=v\vert_{t=0}=w_0.$ Suppose $u \in L^s([0,T];L^q(\mathbb{R}^3))$, $\frac3q +\frac2s =1,$ $q, s \in [2,\infty].$ If $(q,s)=(3,\infty)$, then  $\|u\|_{ L^{\infty}([0,T]; L_{\sigma}^3(\mathbb{R}^3))}$ is assumed sufficiently small. Then  $u\equiv v$.
\end{theorem}

We give the following proposition for which the detailed proof can be seen in Section \ref{w-s}.

\begin{proposition}\label{prop-s-w-uniqueness}
For $p\geq 3,$ $T>0,$ let  ${c}_p$ be as in Theorem \ref{p>3 result}, $|c|>{c}_p$. For $w_0\in L_{\sigma}^p(\mathbb{R}^3)\cap L_{\sigma}^2(\mathbb{R}^3)$, let $w$ be the  $L^p$ mild solution of system (\ref{PNS}) on $[0, T]$. Then $w$ is a $L^2$-weak solution of system (\ref{PNS}) on $[0,T].$
\end{proposition}

According to (\ref{energy-inequality}), there exists $t_0>0$ such that  $w(t_0)\in L_{\sigma}^p\cap L_{\sigma}^3(\mathbb{R}^3)$, $3< p\leq 6$ and $\|w(t_0)\|_{L^3}< \varepsilon_0$. According  to Theorem \ref{p>3 result}, when $|c|>c_p$, there exists a  unique $L^p$ mild solution on $[t_0,\infty)$ to system (\ref{PNS}) with initial data $w(t_0)$.

\begin{corollary}\label{corollary}
For ${w}_0\in L_{\sigma}^2(\mathbb{R}^3)$, let $w$ be a $L^2$-weak solution of system (\ref{PNS}). Then for every $3\leq p\leq 6$ and $|c|>{c}_p$, there exists $T>0$ such that $w(\cdot+T)$ is a $L^p$ mild solution to system (\ref{PNS}) with initial data $w(T)\in L^p_{\sigma}\cap  L_{\sigma}^2(\mathbb{R}^3)$.
\end{corollary}
\begin{remark}\label{remark-Kar}
  Under the condition of Corollary \ref{corollary}, we have $\nabla (|{w}|^{\frac p2})\in L^2 ( [T,\infty);L^2(\mathbb{R}^3))$, and
\begin{equation}
\lim_{t\rightarrow \infty}\|{w}(t)\|_{L^2(\mathbb{R}^3)}=0.
\end{equation}
Furthermore, for $q\geq 3,$ $|c|>\tilde{c}_q$, where $\tilde{c}_q$ is as in Theorem \ref{L^3 well-posedness},
\begin{equation}
\|{w}(t)\|_{L^q(\mathbb{R}^3)} \leq (\frac12-\frac1q)^{\frac32 (\frac12-\frac1q)} (t-T)^{\frac{3}{2}({\frac1q-\frac{1}{2}})}\|{w}(T)\|_{L^2(\mathbb{R}^3)},\quad \text{for all} \quad  t>T.
\end{equation}
\end{remark}
For general initial data $w_0\in L_{\sigma}^3(\mathbb{R}^3)$, we will give the global existence of  $L^2+L^3$ weak solution to system (\ref{PNS}). Inspired by  method in \cite{Cal,Ka,Se}, for any $w_0\in L^3(\mathbb{R}^3)$, we make a  decomposition
\begin{equation}
w_0=v_{10}+v_{20},
\end{equation}
with $\|v_{10}\|_{L^3}<\varepsilon_0$, where $\varepsilon_0$ is as in Theorem \ref{L^3 well-posedness}  and $v_{20}\in L^2\cap L^3(\mathbb{R}^3).$ According to Theorem \ref{L^3 well-posedness}, there exists a unique global $L^3$ mild solution $v_1$ to system
\begin{equation}\label{v1}
\begin{cases}
 \partial_t v_{1}-\Delta v_{1}+(v_{1} \cdot \nabla) v_{1}+(v_{1} \cdot \nabla) v_c+\left(v_c \cdot \nabla\right) v_{1}+\nabla \pi_1=0,\\
 \nabla\cdot v_{1}=0,\\
 v_{1}(x, 0) =v_{10}.
\end{cases}
\end{equation}
Set
\begin{equation}
v_2=w-v_1.
\end{equation}
 Then $v_2$ satisfies
\begin{equation}\label{v2}
\begin{cases}
 \partial_t v_{2}-\Delta v_{2}+(v_{2} \cdot \nabla) v_{2}+(v_{2} \cdot \nabla) (v_c+v_1)+\left((v_c+v_1) \cdot \nabla\right) v_{2}+\nabla \pi_2=0,\\
 \nabla\cdot v_{2}=0,\\
\end{cases}
\end{equation}with $v_{2}(x, 0) =v_{20}$.

We can get the global existence of $w$ by investigating the global existence of $v_2.$ From Theorem 2.7 in \cite{Ka}, system (\ref{v2}) has a weak solution
\begin{equation}\label{v_2-1}
v_2\in C_{w}\left([0, T]; L_{\sigma}^{2}\left(\mathbb{R}^{3}\right)\right) \cap L^{2}\left([0, T]; \dot{H}_{\sigma}^{1}\left(\mathbb{R}^{3}\right)\right)
\end{equation}
for each $T>0$ satisfying the strong energy inequality
\begin{equation}
\|v_2(t)\|_{2}^{2}+2\left(1-K \sup _{t>0}\|v_c+v_1\|_{L^{3}_{w}}\right) \int_{s}^{t}\|\nabla v_2(\tau)\|_{2}^{2} \mathrm{d} \tau \leqslant\|v_2(s)\|_{2}^{2}
\end{equation}
for a constant $K>0$, almost all $s\geq 0$ and all $t\geq s$, and
\begin{equation}\label{v_2-3}
\lim_{t\rightarrow \infty}\|v_2(t)\|_{2}=0.
\end{equation}

In the spirit of the notion of weak $L^3$-solution introduced in Seregin and $\check {\mathrm S}$ver$\acute {\mathrm a}$k \cite{Se}, we give the following definition of $L^2+L^3$ weak solution of system (\ref{PNS}).

\begin{definition}\label{global 3}
Let $T>0$ and $w_0\in L_{\sigma}^3(\mathbb{R}^3)$. A function $w$ is called a $L^2+L^3$ weak solution to system (\ref{PNS}) in $\mathbb{R}^3\times (0,T)$ if $w=v_1+v_2$ for some  $v_1\in C((0,T);L_\sigma^3(\mathbb{R}^3))\cap L^4((0,T);L_\sigma^6(\mathbb{R}^3)),$ and  $v_2\in C_{w}\left((0, T); L_{\sigma}^{2}\left(\mathbb{R}^{3}\right)\right) \cap L^{2}\left((0, T); \dot{H}_{\sigma}^{1}\left(\mathbb{R}^{3}\right)\right)$ such that $v_1$ is a $L^3$ mild solution of (\ref{v1}) and $v_2$ satisfies  the following conditions:\\
$i)$ $v_2$ satisfies  (\ref{v2})
in the sense of distributions; \\
$ii)$
\begin{equation}
w_0=v_1(\cdot, 0)+v_2(\cdot, 0),
\end{equation}
and
\begin{equation}
\lim_{t\rightarrow 0}\|v_2(\cdot, t)-v_2(\cdot, 0)\|_{L^2}=0;
\end{equation}
$iii)$ For all $t\in (0,T)$
\begin{eqnarray}
&&\frac12 \int_{\mathbb{R}^3}|v_2(x,t)|^2dx+\int_0^t \int_{\mathbb{R}^3}|\nabla v_2|^2(x,s)dxds\\
&\leq& \frac 12 \int_{\mathbb{R}^3}|v_{20}(x)|^2dx+\int_0^t \int_{\mathbb{R}^3} v_2\otimes (v_c+v_1): \nabla v_2dxds;\nonumber
\end{eqnarray}
$iv)$ For a.a. $t\in (0,T)$ and any non-negative function $\varphi \in C_{c}^{\infty}\left(\mathbb{R}^3 \times (0,T)\right)$
\begin{eqnarray}
&&\int_{\mathbb{R}^{3}} \left|v_{2}(x, t)\right|^{2} \varphi(x, t)d x+2 \int_{0}^{t} \int_{\mathbb{R}^{3}} \left|\nabla v_{2}\right|^{2} \varphi d x d s \\
&\leq& \int_{0}^{t} \int_{\mathbb{R}^{3}} (2(v_c+v_1)\otimes v_2:\nabla v_2 \varphi + \left(\partial_{s} \varphi + \Delta \varphi\right)|v_2|^{2}\nonumber\\
&&  +(|v_2|^{2}+2 \pi_2+2(v_c+v_1)\cdot v_2)(v_2 \cdot \nabla) \varphi + |v_2|^2 (v_c+v_1)\cdot \nabla \varphi) d x d s.\nonumber
\end{eqnarray}

\end{definition}
We say $w$ is a global $L^2+L^3$ weak solution to system (\ref{PNS}) if it is a $L^2+L^3$ weak solution to system (\ref{PNS}) in $\mathbb{R}^3\times (0,T)$ for all $0<T<\infty.$
Hence, we give the  existence of global $L^2+L^3$ weak solutions to system (\ref{PNS}) as follows
\begin{theorem}  \label{w-global-weak-existence}
Assume that $w_0\in L_{\sigma}^3(\mathbb{R}^3)$ has a decomposition $w_0=v_{10}+v_{20}$ with $v_{10}\in L_{\sigma}^3(\mathbb{R}^3)$, $\|v_{10}\|_{L^3(\mathbb{R}^3)}<\varepsilon_0$ and $v_{20}\in L_{\sigma}^2\cap L_{\sigma}^3(\mathbb{R}^3)$ where  $\varepsilon_0$ is as in Theorem \ref{L^3 well-posedness}.  Then, there exists a global $L^2+L^3$ weak solution $w$ to system (\ref{PNS}) with $w=v_1+v_2$, $v_{1}(\cdot, 0) =v_{10}$ and $v_{2}(\cdot, 0) =v_{20}$.
\end{theorem}
Proof of  Theorem \ref{w-global-weak-existence} is based on   proof of Theorem 2.1 in \cite{Ka}. For the convenience of the reader, we will give details in Section \ref{Sec3}.

\begin{remark}
When initial data $w_0\in L_{\sigma}^p(\mathbb{R}^3)$, $2<p\leq 3,$ by interpolation theory, $w_0$ has a decomposition $w_0=v_{10}+v_{20}$ with $v_{10}\in L_{\sigma}^3(\mathbb{R}^3)$, $\|v_{10}\|_{L^3(\mathbb{R}^3)}<\varepsilon_0$ and $v_{20}\in L_{\sigma}^2$ where  $\varepsilon_0$ is as in Theorem \ref{L^3 well-posedness}. Then, we can easily obtain the global existence of $L^2+L^3$ weak solution to system (\ref{PNS}).

\end{remark}

\textbf{Scheme of the proof and organization of the paper.}
In Section \ref{Sec2}, we give the proof of Theorem \ref{L^3 well-posedness}. In other words, we prove the local well-posedness of solutions to system (\ref{PNS}) with general initial data and global well-posedness of solutions to system (\ref{PNS}) with  small initial data in $L_{\sigma}^3$ space. Also, we investigate the decay rate of solutions to system (\ref{PNS}).
In Section \ref{linear}, properties of linear operator $\mathcal{L}$ on $L^p$, $1<p<\infty$ spaces are studied.
 In Section \ref{w-s}, we prove the Theorem \ref{s-w-uniqueness}, Proposition \ref{prop-s-w-uniqueness} and illustrate Corollary \ref{corollary} briefly. In Section \ref{Sec3}, we illustrate  Theorem \ref{w-global-weak-existence}, i.e. the global existence of $L^2+L^3$ weak solution to system (\ref{PNS}). In Section \ref{Sec4}, we give the proof of Theorem \ref{p>3 result}. In Section \ref{proof of {continuty} }, we give the detailed proof of Theorem \ref{continuty}.

Let us complete this section by the notations that we shall use in this article.

\textbf{Notations.}\\

$\bullet$ We denote $\|\cdot\|_{p}$ or $\|\cdot\|_{L^p}$ the norm of the Lebesgue space $L_x^p(\mathbb{R}^3)$ with $p\in [1,\infty].$

$\bullet$ We denote $\|\cdot\|_{L_t^pL_x^q}$ the norm of the Lebesgue space  $L_t^p([0,\infty); L_x^q(\mathbb{R}^3))$ with $p, q\in [1,\infty].$

$\bullet$ We denote $\|\cdot\|_{C_tL_x^q}$ the norm of the Lebesgue space  $C([0,\infty); L_x^q(\mathbb{R}^3))$ with $q\in [1,\infty].$

$\bullet$ We denote $\|\cdot\|_{L_T^pL_x^q}$ the norm of the Lebesgue space  $L_t^p([0,T]; L_x^q(\mathbb{R}^3))$ with $p, q\in [1,\infty].$

$\bullet$ We denote $\|\cdot\|_{C_TL_x^q}$ the norm of the Lebesgue space  $C_t([0,T]; L_x^q(\mathbb{R}^3))$ with $q\in [1,\infty].$

$\bullet$ $C_0^{\infty}(\mathbb{R}^3)$ denotes  the set of smooth and compactly supported functions.

$\bullet$    $C_w([0,T]; L_x^q(\mathbb{R}^3))$ with $q\in [1,\infty)$ denotes the set of weakly continuous $L^{q}(\mathbb{R}^3)$-valued functions in $t$, i.e. for any $t_{0} \in [0,T]$ and $w \in L^{q^{\prime}}(\mathbb{R}^3)$,
\begin{equation*}
\int_{\mathbb{R}^3} v(x, t) \cdot w(x) d x \rightarrow \int_{\mathbb{R}^3} v\left(x, t_{0}\right) \cdot w(x) d x \quad \text { as } t \rightarrow t_{0}.
\end{equation*}

$\bullet$ For each space $Y,$ we set $Y_{\sigma}=\left\{u \in Y: \text { div } u=0\right\}.$

$\bullet$ We denote $u_i$ the $i$th coordinate ($i=1,2,3$) of a vector $u$.

$\bullet$ Constants independent of solutions may change from line to line and will be denoted by $C$.

\section{Proof of Theorem \ref{L^3 well-posedness}}\label{Sec2}

In this section, we will give the proof of Theorem \ref{L^3 well-posedness}. Our method is based on the following contraction mapping theorem (cf. \cite{Ba}, Theorem 1.72):
\begin{lemma} \label{lem1}
Let $E$ be a Banach space, $N$ be a continuous bilinear map from $E\times E$ to $E$, and $\alpha$ be a positive real number such that
\begin{equation}
\alpha < \frac{1}{4\|N\|} \text{ with } \|N\|:= \sup_{\|u\|,\|v\|\leq 1}\|N(u,v)\|.
\end{equation}
Then for any $a$ in a ball $B(0,\alpha)$ (i.e., with center 0 and radius $\alpha$) in $E$, there exists a unique $x$  in ball $B(0,2\alpha)$ such that
\begin{equation}
x=a+N(x,x).
\end{equation}
\end{lemma}

We will also use a property of Landau solutions $v_c$ which can be obtained by a direct calculation.
\begin{lemma}\label{|x|v_c}
Let $v_c$ be the Landau solutions given by  (\ref{v_c p_c}), then we have
\begin{equation}
\||x|v_c\|_{L^{\infty}}\leq \frac{2\sqrt{2}}{|c|-1}:=K_c.
\end{equation}
\end{lemma}
Next lemma is a fundamental inequality with the singular weight in Sobolev spaces: the so-called Hardy inequality which goes back to the pioneering work by G.H. Hardy in \cite{Har1,Har2}.
\begin{lemma}\label{hardy}
For any $f$ in $\dot{H}^1({\mathbb{R}^{3}}),$ there holds
\begin{equation}
\left(\int_{\mathbb{R}^{3}}\frac{|f(x)|^2}{|x|^2}dx\right)^{\frac{1}{2}}\leq 2\|\nabla f\|_{L^2({\mathbb{R}^{3}})}.
\end{equation}

\end{lemma}

To complete the proof of Theorem \ref{L^3 well-posedness}, we need Lemmas \ref{alemm} and \ref{lem2} which give the results for linear part $a$ and nonlinear part $N$ in (\ref{eq}), separately. Linear part $a$ satisfies the following Cauchy problem
\begin{equation}\label{a}
\begin{cases}
 a_{t}-\Delta a+(a \cdot \nabla) v_{c}+\left(v_{c} \cdot \nabla\right) a+\nabla \pi_1=0,\\
 \nabla\cdot a=0,\\
 a(x, 0) =w_{0}(x).
\end{cases}
\end{equation}
Namely, $a(x,t)$ satisfies
\begin{equation}\label{initial}
\int_{\mathbb{R}^{3}} w_0 \phi d x+ \int_{0}^{\infty} \int_{\mathbb{R}^{3}}\big\{ w\left(-\partial_{t} \phi - \Delta \phi\right)-w\otimes w :\nabla \phi
-(w\otimes v_c+v_c\otimes w):\nabla \phi \big\}dxdt=0,
\end{equation}
for all $\phi\in C_c^{\infty}([0,\infty)\times \mathbb{R}^{3})$ satisfying $\nabla \cdot \phi=0.$

\begin{lemma} \label{alemm}
 For every $c$ satisfies  (\ref{|c|2}), there exists a unique global-in-time  solution $a(x,t)\in C([0,\infty),$ $L_{\sigma}^3(\mathbb{R}^{3}))\cap L^4([0,\infty),L^6_{\sigma}(\mathbb{R}^{3}))$ to system (\ref{a}) with initial data $w_0\in L_{\sigma}^3({\mathbb{R}^{3}})$. Moreover,
\begin{equation}\label{2.7}
\|a(\cdot, t)\|_{L^3}\leq \|a(\cdot, s)\|_{L^3},
\end{equation}
for any $0\leq t\leq s<\infty,$ and
\begin{equation}\label{a-estimate}
\|a\|_{C_tL_x^3\cap L^4_tL^6_x}+\left\|\nabla\left(|a|^{\frac{3}{2}}\right)\right\|^{\frac23}_{L_t^2L_x^2}\leq C_1\|w_0\|_{L^3},
\end{equation}
for a universal constant $C_1.$
\end{lemma}
\noindent\textbf{Proof.}
By classical approximation method, it is easy to obtain the global existence of solutions $a$. For simplicity, we omit the detailed proof and give $a$-$prioi$ estimate for $a.$ Suppose $a$ is sufficiently smooth,
multiplying the equation ($\ref{a})_1$ by $|a|a$, then integrating it on $\mathbb{R}^3$, we have
\begin{eqnarray}\label{prop L^3}
&&\frac{1}{3}\frac{\mathrm{d}}{\mathrm{d} t}\|a(t)\|^3_{L^3}+\frac{8}{9}\|\nabla(|a|^{\frac{3}{2}})\|^2_{L^2}\\\nonumber
&=&-\int_{\mathbb{R}^3}\text{div}(a\otimes v_c+v_c\otimes a)\cdot (|a|a) dx-\int_{\mathbb{R}^3} \nabla \pi_1 \cdot (|a|a) dx.
\end{eqnarray}
For the first term on the right hand side of (\ref{prop L^3}), by using  integration by parts, H\"{o}lder's inequality,  Lemma \ref{|x|v_c} and Hardy inequality in Lemma \ref{hardy}, we have
\begin{eqnarray}\label{prop- rh1}
 -\int_{\mathbb{R}^3}\text{div}(a\otimes v_c+v_c\otimes a)\cdot (|a|a) dx
&=&\int_{\mathbb{R}^3}(a\otimes v_c+v_c\otimes a) \cdot \nabla(|a|a )dx\nonumber\\
&\leq& 4\int_{\mathbb{R}^3}|\nabla a| |a|^{2} |v_c|dx\nonumber\\
&\leq& \frac{8}{3}\int_{\mathbb{R}^3}|\nabla(|a|^{\frac{3}{2}} )||a|^{\frac{3}{2}} |v_c|dx\nonumber\\
&\leq& \frac{8}{3}\left\||x| v_c\right\|_{L^{\infty}}\left\|\nabla(|a|^{\frac{3}{2}} )\right\|_{L^2}\left\|\frac{|a|^{\frac{3}{2}}}{|x|}\right\|_{L^2}\nonumber\\
&\leq& \frac{16}{3}K_c\left\|\nabla(|a|^{\frac{3}{2}} )\right\|_{L^2}^2.
\end{eqnarray}
Then we will estimate the second term on the right hand side of (\ref{prop L^3}).  Thanks to system (\ref{a}), we have
\begin{equation}
\pi_1=-\Delta^{-1} \partial_i\partial_j  \left(a\otimes v_c+v_c\otimes a\right).
\end{equation}
Operator $\Delta^{-1} \partial_i\partial_j $ is Calder{\'o}n-Zygmund operator. According to Example 9.1.7 in \cite{Gra}, there holds $|x|^{p-2}\in A_p$ with $1<p<\infty$. By  H\"{o}lder's inequality, Hardy inequality in Lemma \ref{hardy}, Sobolev embedding and
 boundedness of the Riesz transforms on weighted $L^p$ spaces (see Theorem 9.4.6 in \cite{Gra}), we have
\begin{eqnarray}\label{prop- rh2}
 \int_{\mathbb{R}^3} \nabla \pi_1 \cdot (|a|a) dx
&\leq& \frac23\int_{\mathbb{R}^3} |x|^{\frac13}\left|\Delta^{-1} \partial_i\partial_j  \left(a\otimes v_c+v_c\otimes a\right) \right| \left|\nabla (|a|^{\frac32})\right|\frac{|a|^{\frac12}}{|x|^{\frac13}}dx\nonumber\\
&\leq& \frac23 C_3\||x|^{\frac13}\left(a\otimes v_c+v_c\otimes a\right)\|_{L^{3}} \|\nabla (|a|^{\frac32})\|_{L^{2}} \left\|\frac{|a|^{\frac12}}{|x|^{\frac13}}\right\|_{L^{6}}\nonumber\\
&\leq&  \frac43 C_3\||x|v_c\|_{L^{\infty}} \left\|\frac{a}{|x|^{\frac23}}\right\|_{L^3}\|\nabla (|a|^{\frac32})\|_{L^{2}} \left\|\nabla(|a|^{\frac32})\right\|_{L^{2}}^{\frac13}\nonumber\\
&\leq& \frac83 C_3 K_c\|\nabla (|a|^{\frac{3}{2}})\|_{L^2}^2.
\end{eqnarray}
Combining (\ref{prop L^3}), (\ref{prop- rh1}) and (\ref{prop- rh2}), we deduce
\begin{equation}\label{2.13}
\frac{\mathrm{d}}{\mathrm{d} t}\|a(t)\|^3_{L^3}+\left(\frac{8}{3}-16 K_c-8C_3K_c\right)\left\|\nabla\left(|a|^{\frac{3}{2}}\right)\right\|^2_{L^2}\leq 0.
\end{equation}
Choosing $|c|$ big enough such that
\begin{equation}\label{|c|2}
\frac83-16 K_c-8C_3K_c>0,
\end{equation}
 then we have
\begin{equation}
\sup_t\|a(t)\|^3_{L^3}+\left(\frac{8}{3}-16 K_c-8C_3K_c\right)\left\|\nabla\left(|a|^{\frac{3}{2}}\right)\right\|^2_{L_t^2L_x^2}\leq \|w_0\|^3_{L^3}.
\end{equation}
Hence, there exists a constant $C$ such that
\begin{equation}
\sup_t\|a(t)\|_{L^3}+\left\|\nabla\left(|a|^{\frac{3}{2}}\right)\right\|^{\frac23}_{L_t^2L_x^2}\leq C\|w_0\|_{L^3}.
\end{equation}
Hence, we deduce (\ref{2.7}).
By interpolation theory, we have
\begin{equation}
 \|a\|_{L_t^4L_x^6}\leq \|a\|^{\frac14}_{L_t^{\infty}L_x^3}\left\|\nabla\left(|a|^{\frac{3}{2}}\right)\right\|^{\frac12}_{L_t^2L_x^2}
\leq C\left(\sup_t\|a(t)\|_{L^3}+\left\|\nabla(|a|^{\frac{3}{2}})\right\|_{L_t^2L_x^2}^{\frac{2}{3}}\right)
 \leq C\|w_0\|_{L^3}.
\end{equation}
Therefore, there holds
\begin{equation}\label{a-prioi}
\|a\|_{L_t^\infty L_x^3\cap L^4_tL^6_x} +\left\|\nabla\left(|a|^{\frac{3}{2}}\right)\right\|^{\frac23}_{L_t^2L_x^2}\leq C_1\|w_0\|_{L^3},
\end{equation}
for a constant $C_1$.

Then, we consider the  continuity of  solution $a$ over time $t$. Because of the translational invariance in time, we only consider time around 0. For any sequence $t_{n}\rightarrow 0,$ according to   (\ref{a-prioi}),
there holds
 \begin{equation}\label{weak con}
a\in L_t^\infty L_x^3.
 \end{equation}
 Therefore, there exist subsequence $t_{n_j}\rightarrow 0$ such that
\begin{equation*}
a(\cdot,t_{n_j})\rightharpoonup w_0(\cdot)  \text{\ weakly in } L^3.
  \end{equation*}
  Therefore, we have
 \begin{equation}\label{1}
\|w_0\|_{L^3}\leq    {\underline{\lim}}_{j \rightarrow \infty}\|a(\cdot,t_{n_j} )\|_{L^3}.
 \end{equation}
By energy inequality (\ref{2.13}), there holds
\begin{equation}\label{20}
  {\overline{\lim}}_{j \rightarrow\infty}\|a(\cdot,t_{n_j} )\|_{L^3}\leq \|w_0\|_{L^3} .
 \end{equation}
Combining with (\ref{1}) and (\ref{20}), there holds
\begin{equation} \label{con_1}
\lim_{j \rightarrow\infty}\|a(\cdot,t_{n_j} )\|_{L^3}=\|w_0\|_{L^3}.
\end{equation}
Hence
\begin{equation}
a(\cdot, t_{n_j})\rightarrow w_0  \text{\  in } L^3.
  \end{equation}
  Then, we deduce
  \begin{equation}
a(\cdot,t)\rightarrow w_0(\cdot)  \text{\  in } L^3 \text{\  as }t\rightarrow 0.
  \end{equation}
Hence, $a\in C([0,\infty);L_x^3)$.
 {\hfill
$\square$\medskip}
\begin{remark}
Indeed, more strictly, we can  prove the  existence of   $a$  by approximation theory.
Assume $a_0=0$, we construct the iterative sequence $\{a_k\} $ as follows
\begin{equation*}\label{a-k}
\begin{cases}
&\partial_t a_{k}-\Delta a_{k}=-(a_{k-1} \cdot \nabla) v_{c}-\left(v_{c} \cdot \nabla\right) a_{k-1}-\nabla \pi_{k-1},\ \ \text{for\ \ }k=1,2,\cdots,\\
&\nabla \cdot a_k=0, \\
&\pi_{k-1}= (-\Delta)^{-1}\partial_i \partial_j(v_c\otimes a_{k-1}+a_{k-1}\otimes v_c), 
\\ &a_{k}|_{t=0}=w_0.
 \end{cases}
 \end{equation*}
 We claim that $\{a_{k}\}\in C([0,\infty);L_x^3)\cap L_t^4([0,\infty);L_x^6)$ and $\nabla(|a_k|^{\frac32})\in L_t^2([0,\infty);L_x^2).$ By Duhamel principle, $a_k$ also satisfy the integral formulation
$e^{t\Delta}a_k-\int_0^t e^{(t-s)\Delta}\mathbb{P}\text{div}(a_{k-1} \otimes v_c+v_c\otimes a_{k-1} )ds.$ Since semigroup $e^{t\Delta}: L^3\rightarrow L^3,$ we have $a_k\in  C([0,\infty);L_x^3).$  By energy estimate, we have $\{a_k\}$ is a Cauchy sequence in $L_t^4([0,\infty);L_x^6)$. Limit of $\{a_k\}$ is $a$ which satisfies Lemma \ref{alemm}. We omit the details.
\end{remark}

For $w_0\in L_{\sigma}^3(\mathbb{R}^3)$ and  $0\leq t<\infty,$ let
\begin{equation}
T(t)w_0:=a(x,t),
\end{equation}
where $a(x,t)$ is the unique solution of (\ref{a}) given by Lemma \ref{alemm}. Then $T(t),$ $0\leq t<\infty,$ is a one parameter family of bounded linear operators from $L_{\sigma}^3(\mathbb{R}^3)$ into $L_{\sigma}^3(\mathbb{R}^3)$ satisfying $T(0)=I,$ the identity operator of $L_{\sigma}^3(\mathbb{R}^3)$, $T(t+s)=T(t)T(s)$ for every $t,s \geq 0,$ $\|T(t)\|\leq 1$ for every $t\geq 0,$ and $\lim_{t\rightarrow 0+}T(t)w_0=w_0$ in $L_{\sigma}^3(\mathbb{R}^3)$.

Therefore, $T(t)$  is a strongly continuous semigroup of contraction, see Definition 1.2.1 and Section 1.3 in \cite{Paz}.

The linearized operator $-\mathcal{L}$, given in (\ref{L def}), with the domain of definition
\begin{equation}
D(-\mathcal{L}):= \Big\{w_0\in L_{\sigma}^3(\mathbb{R}^3): \lim_{t\rightarrow 0+}\frac{T(t)w_0-w_0}{t} \text{ exists in }L_{\sigma}^3(\mathbb{R}^3)\Big\},
\end{equation}
is the infinitesimal generator of the  semigroup $T(t),$ see Section 1.1 of \cite{Paz}. By Corollary 1.2.5 in \cite{Paz}, $D(-\mathcal{L})$ is dense in $L_{\sigma}^3(\mathbb{R}^3)$ and
$-\mathcal{L}$ is a closed linear operator in $L_{\sigma}^3(\mathbb{R}^3)$. We also denote $T(t)$ as $e^{-t\mathcal{L}}$.

 Next, we will estimate the nonlinear part $N(w_1,w_2)$. Denote $z=N(w_1,w_2),$ it's obvious that $z$ satisfies the following system
\begin{eqnarray}\label{N}
\begin{cases}{}
 z_{t}-\Delta z+(z \cdot \nabla) v_{c}+\left(v_{c} \cdot \nabla\right) z+\nabla \pi_2=-\text{div}(w_1\otimes w_2),\\
 \nabla\cdot z=0,\\
 z(x, 0) =0.
\end{cases}
\end{eqnarray}
\begin{lemma}\label{lem2}
 For every $c$ satisfies (\ref{|c|3}),
there exists a   unique   solution $z(x,t)\in C([0,T],L_{\sigma}^3(\mathbb{R}^{3}))\cap L^4([0,T],L^6_{\sigma}(\mathbb{R}^{3}))$ to  system (\ref{N}) with $w_1, w_2\in L^4([0,T];L^6(\mathbb{R}^3))$, satisfying
\begin{equation}\label{z-estimate}
 \|z(t)\|_{C([0,T];L_x^3)\cap L^4_t([0,T];L^6_x)}+\left\|\nabla\left(|z|^{\frac{3}{2}}\right)\right\|^{\frac23}_{L_t^2([0,T];L_x^2)}
\leq C_2 \|w_1\|_{L_t^4([0,T];L_x^6)}\|w_2\|_{L_t^4([0,T];L_x^6)},
\end{equation}
for a universal constant $C_2$ which is independent of  $T$.
\end{lemma}
\noindent\textbf{Proof.}
We omit the detailed proof of the existence of solution $z$ since it can be obtained by classical approximation method. Then, we give $a$-$prioi$ estimate for $z.$ Suppose $z$ is sufficiently smooth,
multiplying the equation $(\ref{N})_1$ by $|z|z$ and integrating it on $\mathbb{R}^3$, we have
\begin{eqnarray*}
&&\frac{1}{3}\frac{\mathrm{d}}{\mathrm{d} t}\|z(t)\|^3_{L^3}+\frac{8}{9}\|\nabla(|z|^{\frac{3}{2}})\|^2_{L^2}\\
&=&-\int_{\mathbb{R}^3}\text{div}(z\otimes v_c+v_c\otimes z)\cdot (|z|z)dx-\int_{\mathbb{R}^3}\text{div}(w_1\otimes w_2)\cdot (|z|z)dx
-\int_{\mathbb{R}^3} \nabla \pi_2 \cdot (|z|z) dx.\nonumber
\end{eqnarray*}
The estimate for the first term on the right hand side  is the same as (\ref{prop- rh1}). Hence, there holds
\begin{equation}\label{aaaaa}
\begin{aligned}
-\int_{\mathbb{R}^3}\text{div}(z\otimes v_c+v_c\otimes z)\cdot  (|z|z) dx\leq \frac{16}{3}K_c\|\nabla(|z|^{\frac{3}{2}} )\|_{L^2}^2.\\
\end{aligned}
\end{equation}
For the second term, by integration by parts,  H\"{o}lder's inequality and Young's inequality, we have
\begin{eqnarray}
-\int_{\mathbb{R}^3}\text{div}(w_1\otimes w_2) \cdot (|z|z) dx &=&\int_{\mathbb{R}^3}(w_1\otimes w_2)\cdot  \nabla(|z|z )dx\nonumber\\
&\leq& \frac43\int_{\mathbb{R}^3}\left|w_1\otimes w_2\right|
\left|\nabla\left(|z|^{\frac{3}{2}} \right)\right||z|^{\frac{1}{2}} dx\nonumber\\
&\leq& \frac43\left\|\nabla\left(|z|^{\frac{3}{2}} \right)\right\|_{L^2} \left\|(w_1\otimes w_2)|z|^{\frac{1}{2}}\right\|_{L^2}\nonumber\\
&\leq& \frac43\left\|\nabla\left(|z|^{\frac{3}{2}} \right)\right\|_{L^2}\left\||z|^{\frac{1}{2}}\right\|_{L^6}\left\|w_1\otimes w_2\right\|_{L^3}\nonumber\\
&\leq& \frac43\left\|\nabla\left(|z|^{\frac{3}{2}} \right)\right\|_{L^2}\left\|z\right\|_{L^3}^{\frac{1}{2}}\left\|w_1\otimes w_2\right\|_{L^3}\\
&\leq& \frac{2}{15}\left\|\nabla\left(|z|^{\frac{3}{2}} \right)\right\|^2_{L^2}+\frac{10}{3}\left\|z\right\|_{L^3}\left\|w_1\otimes w_2\right\|^2_{L^3}.\nonumber
\end{eqnarray}
Since $w_1, w_2\in  L_T^{4} L_x^6$, we have
\begin{eqnarray}\label{prop rh2_3}
&&-\int_0^T\int_{\mathbb{R}^3}\text{div}(w_1\otimes w_2) \cdot (|z|z) dxdt \nonumber\\
&\leq& \frac{2}{15}\int_0^T\left\|\nabla\left(|z|^{\frac{3}{2}} \right)\right\|^2_{L^2}dt+\frac{10}{3} \left\|z\right\|_{L_T^{\infty}L_x^3}\left\|w_1\right\|^2_{L_T^4L_x^6}\left\|w_2\right\|^2_{L_T^4L_x^6}.
\end{eqnarray}
For the third term on the right hand,  according to system (\ref{N}), we have
\begin{equation}
\pi_2=-\Delta^{-1} \partial_i\partial_j  \left(z\otimes v_c+v_c\otimes z+w_1 \otimes w_2\right).
\end{equation}
Using  integration by parts and H\"{o}lder's inequality, we have the following estimate
\begin{eqnarray}\label{pi}
&&\int_{\mathbb{R}^3} \nabla \pi_2 \cdot (|z|z) dx
=-\int_{\mathbb{R}^3} \pi_2 \nabla(|z|)\cdot z dx\nonumber\\
&=&\int_{\mathbb{R}^3} \Delta^{-1} \partial_i\partial_j  \left(z\otimes v_c+v_c\otimes z+w_1 \otimes w_2\right)\nabla(|z|)\cdot z dx\nonumber\\
&\leq& \frac23 \int_{\mathbb{R}^3} \left|\Delta^{-1} \partial_i\partial_j  \left(z\otimes v_c+v_c\otimes z+w_1 \otimes w_2\right)\right|\left|
\nabla\left(|z|^{\frac{3}{2}} \right)\right||z|^{\frac{1}{2}} dx\nonumber\\
&\leq& \frac23 \int_{\mathbb{R}^3} \left|\Delta^{-1} \partial_i\partial_j  \left(z\otimes v_c+v_c\otimes z\right)\right|\left|
\nabla\left(|z|^{\frac{3}{2}} \right)\right||z|^{\frac{1}{2}} dx\nonumber\\
&&+\frac23 \int_{\mathbb{R}^3} \left|\Delta^{-1} \partial_i\partial_j  \left(w_1 \otimes w_2\right)\right|\left|
\nabla\left(|z|^{\frac{3}{2}} \right)\right||z|^{\frac{1}{2}} dx.
\end{eqnarray}
Thanks to (\ref{prop- rh2}), we have
\begin{equation}\label{pi1}
\frac23 \int_{\mathbb{R}^3} \left|\Delta^{-1} \partial_i\partial_j  \left(z\otimes v_c+v_c\otimes z\right)\right|\left|
\nabla\left(|z|^{\frac{3}{2}} \right)\right||z|^{\frac{1}{2}} dx
\leq \frac{8}{3} C_3 K_c \|\nabla(|z|^{\frac{3}{2}} )\|_{L^2}^{2}.
\end{equation}
For the second part, by H\"{o}lder's inequality, property of Riesz operator and Young's inequality, we have
\begin{eqnarray}\label{pi2}
&&\frac23 \int_{\mathbb{R}^3} \left|\Delta^{-1} \partial_i\partial_j  \left(w_1\otimes w_2\right)\right|\left|\nabla\left(|z|^{\frac{3}{2}} \right)\right||z|^{\frac{1}{2}} dx\nonumber\\
&\leq& \frac23 \| \Delta^{-1} \partial_i\partial_j  \left(w_1\otimes w_2\right)\|_{L^3}\|\nabla(|z|^{\frac{3}{2}} )\|_{L^2}\| |z|^{\frac{1}{2}}\|_{L^6}\nonumber\\
&\leq& \frac23 H_3\| w_1\otimes w_2\|_{L^3}\|\nabla(|z|^{\frac{3}{2}} )\|_{L^2}\| z\|^{\frac12}_{L^3}\nonumber\\
&\leq& \frac{2}{15}H_3\left\|\nabla\left(|z|^{\frac{3}{2}} \right)\right\|^2_{L^2}+\frac{10}{3}H_3\left\|z\right\|_{L^3}\left\|w_1\otimes w_2\right\|^2_{L^3}.
\end{eqnarray}
where $H_3$ is a constant  and origins from the following inequality
\begin{equation}\label{Hr}
\|\Delta^{-1} \partial_i\partial_j f \|_{L^r}\leq H_r \|f\|_{L^r},
\end{equation}
for $1<r<\infty.$ For scalar Riesz transforms, Iwaniec and  Martin \cite{Iwa} showed that the norm $\|R_l\|_{L^r}$ of the Riesz operator $R_l:$ $L^r(\mathbb{R}^n)\rightarrow L^r(\mathbb{R}^n)$ is equal to
\begin{equation}\label{crvalue}
\begin{cases}
\text{tan} (\frac{\pi}{2r}), \text{if   } 1<r\leq2,\\
\text{cot} (\frac{\pi}{2r}), \text{if   } 2\leq r < \infty.
\end{cases}
\end{equation}
Combining with (\ref{pi}), (\ref{pi1}) and (\ref{pi2}), we have the following estimate
\begin{eqnarray*}
\int_{\mathbb{R}^3} \nabla \pi_2 \cdot (|z|z) dx
&\leq& \left(\frac{8}{3} C_3 K_c +\frac{2}{15}H_3\right)\left\|\nabla (|z|^{\frac{3}{2}})\right\|_{L^2}^2+\frac{10}{3}H_3\|z\|_{L^3}|\|w_1\otimes w_2\|^2_{L^3}.
\end{eqnarray*}
Therefore, we deduce
\begin{equation}
\begin{aligned}
& \left\|z\right\|^3_{L_T^{\infty}L_x^3}+\left(\frac83-16K_c-\frac{2}{5}-8C_3 K_c -\frac{2}{5}H_3\right)\left\|\nabla(|z|^{\frac{3}{2}})\right\|^2_{L_T^2L_x^2}\\
&\leq (10+10H_3) \|z\|_{L_T^{\infty}L_x^3}\|w_1\|_{L_T^{4}L_x^6}^2\|w_2\|_{L_T^{4}L_x^6}^2.
\end{aligned}
\end{equation}
Choosing $|c|$ big enough such that
\begin{equation}\label{|c|3}
\frac83-16K_c-\frac{2}{5}-8C_3 K_c -\frac{2}{5}H_3>0,
\end{equation}
we have
\begin{equation}\label{2}
\begin{aligned}
 \left\|z\right\|_{L_T^{\infty}L_x^3}+\left\|\nabla(|z|^{\frac{3}{2}})\right\|^{\frac{2}{3}}_{L_T^2L_x^2}\leq C \|w_1\|_{L_T^{4}L_x^6}\|w_2\|_{L_T^{4}L_x^6},
\end{aligned}
\end{equation}
for a universal constant $C.$
By interpolation theory, there holds
\begin{equation}\label{zzzz}
\|z(t)\|_{L_T^\infty L_x^3\cap L^4_T L^6_x}+\left\|\nabla\left(|z|^{\frac{3}{2}}\right)\right\|^{\frac23}_{L_T^2L_x^2}\leq C_2 \|w_1\|_{L_T^4 L_x^6}\|w_2\|_{L_T^4 L_x^6},
\end{equation}
for a universal constant $C_2.$

By similar argument as (\ref{con_1}), we can obtain the continuity of $z$ over time $t.$ Combining with (\ref{zzzz}), we deduce (\ref{z-estimate}).
 {\hfill
$\square$\medskip}

\noindent\textbf{Proof of Theorem \ref{L^3 well-posedness}.}
For every $c$ satisfies (\ref{|c|2}) and (\ref{|c|3}), according to Lemma  \ref{alemm}, we have
\begin{equation}
\|a(t)\|_{L^4_t L_x^6}\leq C_1 \|w_0\|_{L^3}.
\end{equation}
Applying Lemma  \ref{lem2} with $(w_1,w_2)=(w,w)$, we have
\begin{equation}
\|N\|\leq C_2.
\end{equation}
Using Lemma \ref{lem1} with $E={L^4_tL_x^6}$, when $\|w_0\|_{L^3}< \frac{1}{4C_1C_2}:=\varepsilon_0$, there exists a global unique solution $w\in {L^4_tL_x^6}$ and $\|w\|_{{L^4_t L_x^6}}\leq C \|w_0\|_{L^3}.$ According to Lemmas \ref{alemm} and \ref{lem2}, solution $w\in C([0,\infty);L^3(\mathbb{R}^3))$, $\nabla (|w|^{\frac32})\in L^2([0,\infty);L^2(\mathbb{R}^3))$.

When $w_0\in{L^3(\mathbb{R}^3)}$, thanks to Lemma \ref{alemm},  we have
\begin{equation}
\|a\|_{L^4_tL^6_x} \leq C_1 \|w_0\|_{L^3(\mathbb{R}^3)}.
\end{equation}
There exists $T>0$ such that $\|a\|_{L^4_TL^6_x}< \frac{1}{4C_2}.$ Using Lemma \ref{lem1} and \ref{lem2} with $E={L^4_T L_x^6}$, a unique local solution $w\in L^4_T L_x^6$ exists on $[0,T].$ According to   Lemmas \ref{alemm} and \ref{lem2}, solution $w\in C([0,T];L^3(\mathbb{R}^3))$ and $\nabla (|w|^{\frac32})\in L^2([0,T];L^2(\mathbb{R}^3))$.

Next, we will investigate the decay rate of solution $w$, i.e. (\ref{nonlinear decay}). Our method is inspired by \cite{Ca}. Let $T>0$ and $3 \leq q < \infty$. Denote $r(t)=\frac{1}{\frac{1}{T}(\frac{1}{q}-\frac{1}{3})t+\frac{1}{3}}$. First, we give $a$-$prioi$ estimate of $\|w(\cdot,t)\|_{r(t)}$.
By direct computation, we have
\begin{eqnarray}\label{remk_1}
&&r(t)^{2}\|w(\cdot,t)\|_{r(t)}^{r(t)-1} \frac{\mathrm{d}}{\mathrm{d} t}\|w(\cdot,t)\|_{r(t)}\nonumber\\
&=&\dot{r}(t) \int_{\mathbb{R}^{3}} |w(\cdot,t)|^{r(t)} \ln \left(|w(\cdot,t)|^{r(t)} /\|w(\cdot,t)\|_{r(t)}^{r(t)}\right) \mathrm{d} x\nonumber\\
&&+r(t)^{2} \int_{\mathbb{R}^{3}} |w(\cdot,t)|^{r(t)-1} \frac{\mathrm{d}}{\mathrm{d} t} \left|w(\cdot,t)\right| \mathrm{d}x\nonumber\\
&=&\dot{r}(t) \int_{\mathbb{R}^{3}} |w(\cdot,t)|^{r(t)} \ln \left(|w(\cdot,t)|^{r(t)} /\|w(\cdot,t)\|_{r(t)}^{r(t)}\right) \mathrm{d} x\nonumber\\
&&+r(t)^{2} \int_{\mathbb{R}^{3}} |w(\cdot,t)|^{r(t)-2} w_i\frac{\mathrm{d}}{\mathrm{d} t} w_i  \mathrm{d}x\nonumber\\
&=&\dot{r}(t) \int_{\mathbb{R}^{3}} |w(\cdot,t)|^{r(t)} \ln \left(|w(\cdot,t)|^{r(t)} /\|w(\cdot,t)\|_{r(t)}^{r(t)}\right) \mathrm{d} x\nonumber\\
&& +r(t)^{2} \int_{\mathbb{R}^{3}} |w(\cdot,t)|^{r(t)-2} w_i\left( \partial_j\partial_j w_i-\partial_j (w_iw_j+w_iv_j+v_iw_j)-\partial_i\pi\right)  \mathrm{d}x,\nonumber
\end{eqnarray}
where the last equality holds on account of $(\ref{PNS})_1$.

Set
\begin{equation}
\begin{aligned}
\uppercase\expandafter{\romannumeral1}&=r(t)^{2} \int_{\mathbb{R}^{3}} |w(\cdot,t)|^{r(t)-2}w_i\partial_j\partial_j w_i \mathrm{d}x,\\
\uppercase\expandafter{\romannumeral2}&= -r(t)^{2} \int_{\mathbb{R}^{3}} |w(\cdot,t)|^{r(t)-2} w_i\partial_j (w_iw_j)\mathrm{d}x,\\
\uppercase\expandafter{\romannumeral3}&= -r(t)^{2} \int_{\mathbb{R}^{3}} |w(\cdot,t)|^{r(t)-2} w_i\partial_j (w_iv_j)\mathrm{d}x,\\
\uppercase\expandafter{\romannumeral4}&= -r(t)^{2} \int_{\mathbb{R}^{3}} |w(\cdot,t)|^{r(t)-2} w_i\partial_j (v_iw_j)\mathrm{d}x,\\
\uppercase\expandafter{\romannumeral5}&= -r(t)^{2} \int_{\mathbb{R}^{3}} |w(\cdot,t)|^{r(t)-2}w_i\partial_i\pi \mathrm{d}x.
\end{aligned}
\end{equation}
Using  integration by parts, we have
\begin{eqnarray*}\label{I_1}
\uppercase\expandafter{\romannumeral1}&=& -r(t)^{2} \int_{\mathbb{R}^{3}} \partial_j(|w(t)|^{r(t)-2} w_i) \partial_j w_idx \nonumber\\
&=& -r(t)^{2} \int_{\mathbb{R}^{3}} \frac12\partial_j|w(t)|^{r(t)-2} \partial_j |w|^2dx-r(t)^{2} \int_{\mathbb{R}^{3}}  |w(t)|^{r(t)-2}|\nabla w|^2dx \nonumber\\
&=& -r(t)^{2} \int_{\mathbb{R}^{3}} \frac{4(r(t)-2) }{r(t)^{2} }\left|\nabla (|w(\cdot,t)|^{r(t) / 2})\right|^{2}dx -r(t)^{2} \int_{\mathbb{R}^{3}}  |w(t)|^{r(t)-2}|\nabla w|^2dx\nonumber\\
&=& -4r(t)(r(t)-2)\|\nabla (|w(\cdot,t)|^{\frac{r(t)}{2}})\|_{L^2}^2-r(t)^{2} \int_{\mathbb{R}^{3}}  |w(t)|^{r(t)-2}|\nabla w|^2dx,\nonumber
\end{eqnarray*}
where the third equality holds by use of the fact
\begin{equation*}
\nabla\left(|w(\cdot, t)|^{r(t)-2}\right) \cdot \nabla (|w(\cdot,t)|^2)=\frac{8(r(t)-2) }{r(t)^{2} }\left|\nabla (|w(\cdot,t)|^{r(t) / 2})\right|^{2}.
\end{equation*}
We have
\begin{eqnarray}\label{important_1}
&&r(t)^{2}\|w(\cdot, t)\|_{r(t)}^{r(t)-1} \frac{\mathrm{d}}{\mathrm{d} t}\|w(\cdot, t)\|_{r(t)}\nonumber\\
&=&\dot{r}(t) \int_{\mathbb{R}^{3}} |w(\cdot, t)|^{r(t)} \ln \left(|w(\cdot, t)|^{r(t)} /\|w(\cdot, t)\|_{r(t)}^{r(t)}\right) \mathrm{d}x-4r(t)(r(t)-2)\|\nabla (|w(\cdot,t)|^{\frac{r(t)}{2}})\|_{L^2}^2\nonumber\\
&&-r(t)^{2} \int_{\mathbb{R}^{3}}  |w(\cdot, t)|^{r(t)-2}|\nabla w(\cdot, t)|^2dx +\uppercase\expandafter{\romannumeral2}+\uppercase\expandafter{\romannumeral3}+\uppercase\expandafter{\romannumeral4}+\uppercase\expandafter{\romannumeral5}.
\end{eqnarray}
Next, we will estimate $\uppercase\expandafter{\romannumeral2}-\uppercase\expandafter{\romannumeral5}$, respectively.

Thanks to integration by parts, H\"{o}lder's inequality and Sobolev embedding $\dot{H}^{1}(\mathbb{R}^{3})\hookrightarrow L^6(\mathbb{R}^{3})$ (best constant can be seen in  \cite{Swa}), we have
\begin{eqnarray}
\uppercase\expandafter{\romannumeral2}&=& -\frac{r(t)^{2} }{2}\int_{\mathbb{R}^{3}} |w(\cdot,t)|^{r(t)-2} w_j\partial_j|w|^2\mathrm{d}x\nonumber\\
&=& \frac{r(t)^{2} }{2}\int_{\mathbb{R}^{3}} \partial_j(|w(\cdot,t)|^{r(t)-2}) w_j|w|^2\mathrm{d}x\nonumber\\
&\leq&  r(t)(r(t)-2)\int_{\mathbb{R}^{3}} \left| \nabla(|w(\cdot,t)|^{\frac{r(t)}{2}})\right| |w|^{\frac{r(t)}{2}}|w|\mathrm{d}x\nonumber\\
&\leq&  r(t)(r(t)-2)\left\|\nabla\left(|w(t)|^{\frac{r(t)}{2}}\right)\right\|_{L^{2}}\left\||w(t)|^{\frac{r(t)}{2}}\right\|_{L^{6}}\|w(t)\|_{L^{3}}\nonumber\\
&\leq&  r(t)(r(t)-2)\left\|\nabla\left(|w(t)|^{\frac{r(t)}{2}}\right)\right\|_{L^{2}}^{2}\|w(t)\|_{L^{3}}.
\end{eqnarray}
Combining (\ref{L^3 stability-1}) with $\|w_0\|_{L^3}\leq \varepsilon_0,$ there holds
\begin{equation}\label{I_2}
\begin{aligned}
\uppercase\expandafter{\romannumeral2}\leq r(t)(r(t)-2)C\varepsilon_0 \left\|\nabla (|w(t)|^{\frac{r(t)}{2}})\right\|^2_{L^2}.
\end{aligned}
\end{equation}
According to H\"{o}lder's inequality, the Hardy inequality in Lemma \ref{hardy} and Lemma \ref{|x|v_c}, we deduce
\begin{eqnarray}\label{I_3}
\uppercase\expandafter{\romannumeral3}&=& \frac{r(t)^{2}}{2} \int_{\mathbb{R}^{3}} \partial_j(|w(\cdot,t)|^{r(t)-2} )|w|^2v_j\mathrm{d}x\nonumber\\
&=& r(t)(r(t)-2)\int_{\mathbb{R}^{3}} v_c\cdot \nabla (|w(\cdot,t)|^{\frac{r(t)}{2}})|w(\cdot,t)|^{\frac{r(t)}{2}}\mathrm{d}x\nonumber\\
&\leq&  r(t)(r(t)-2)\left\|\nabla\left(|w(t)|^{\frac{r(t)}{2}}\right)\right\|_{L^{2}}\left\|\frac{|w(t)|^{\frac{r(t)}{2}}}{|x|}\right\|_{L^{2}}\left\||x| v_{c}\right\|_{L^{\infty}}\nonumber\\
&\leq& 2 r(t)(r(t)-2) K_{c}\left\|\nabla\left(|w(t)|^{\frac{r(t)}{2}}\right)\right\|_{L^{2}}^{2}.
\end{eqnarray}
For forth term, we  integrate by parts to have
\begin{eqnarray*}
\uppercase\expandafter{\romannumeral4}&=& -r(t)^{2} \int_{\mathbb{R}^{3}} |w(\cdot,t)|^{r(t)-2} w_i\partial_j (v_iw_j)\mathrm{d}x\\
&=&r(t)^{2}\int_{\mathbb{R}^{3}} \partial_j (|w(\cdot,t)|^{r(t)-2}) w_iv_iw_j\mathrm{d}x
+ r(t)^{2}\int_{\mathbb{R}^{3}} |w(\cdot,t)|^{r(t)-2} \partial_jw_i v_iw_j\mathrm{d}x.
\end{eqnarray*}
The estimate of the first part is similar as (\ref{I_3}). We have
\begin{eqnarray*}
r(t)^{2}\int_{\mathbb{R}^{3}} \partial_j (|w(\cdot,t)|^{r(t)-2}) w_iv_iw_j\mathrm{d}x
\leq 4 r(t)(r(t)-2) K_{c}\left\|\nabla\left(|w(t)|^{\frac{r(t)}{2}}\right)\right\|_{L^{2}}^{2}.
\end{eqnarray*}
By Lemma \ref{|x|v_c}, Cauchy inequality and the Hardy inequality, we can estimate the second part as follows
\begin{eqnarray*}
&&r(t)^{2}\int_{\mathbb{R}^{3}} |w(\cdot,t)|^{r(t)-2} \partial_jw_i v_iw_j\mathrm{d}x\nonumber\\
&\leq& r(t)^{2}\int_{\mathbb{R}^{3}} |w(\cdot,t)|^{r(t)-2} |\partial_jw_i | ||x|v_i|\frac{|w_j|}{|x|}\mathrm{d}x\nonumber\\
&\leq& r(t)^{2}K_c\int_{\mathbb{R}^{3}} |w(\cdot,t)|^{r(t)-2} |\partial_jw_i | \frac{|w_j|}{|x|}\mathrm{d}x\nonumber\\
&\leq& \frac{r(t)^{2}}{2} K_c\int_{\mathbb{R}^{3}} |w(\cdot,t)|^{r(t)-2} |\nabla w(\cdot,t)|^{2} dx+ \frac{r(t)^{2}}{2} K_c
\int_{\mathbb{R}^{3}} |w(\cdot,t)|^{r(t)-2} \frac{|w(\cdot,t)|^{2}}{|x|^2} dx\nonumber\\
&\leq& \frac{r(t)^{2}}{2} K_c\int_{\mathbb{R}^{3}} |w(\cdot,t)|^{r(t)-2} |\nabla w(\cdot,t)|^{2} dx+ 2r(t)^{2} K_c\left\|\nabla\left(|w(t)|^{\frac{r(t)}{2}}\right)\right\|_{L^{2}}^{2}.
\end{eqnarray*}
Therefore, we have
\begin{equation}\label{I_4}
\begin{aligned}
\uppercase\expandafter{\romannumeral4}\leq & \frac{r(t)^{2}}{2} K_c\int_{\mathbb{R}^{3}} |w(\cdot,t)|^{r(t)-2} |\nabla w(\cdot,t)|^{2} dx\\
&+ (4 r(t)(r(t)-2)+2r(t)^{2}) K_c\left\|\nabla\left(|w(t)|^{\frac{r(t)}{2}}\right)\right\|_{L^{2}}^{2}.
\end{aligned}
\end{equation}
Note that the pressure $\pi=-\frac{\partial_i\partial_j}{\Delta}(w_iw_j+v_iw_j+w_iv_j )$, using integration by parts, we obtain
\begin{eqnarray}\label{I_5}
\uppercase\expandafter{\romannumeral5}
&=& r(t)^{2} \int_{\mathbb{R}^{3}} \partial_i(|w(\cdot,t)|^{r(t)-2})w_i\pi \mathrm{d}x\nonumber\\
&=& r(t)^{2} \int_{\mathbb{R}^{3}} \partial_i(|w(\cdot,t)|^{r(t)-2})w_i\left(-\frac{\partial_i\partial_j}{\Delta}(v_iw_j+w_iv_j+w_iw_j)\right) \mathrm{d}x.
\end{eqnarray}
This term is more complex to deal with, we will  estimate it more carefully.

Set
$$
\uppercase\expandafter{\romannumeral5}_1= r(t)^{2} \int_{\mathbb{R}^{3}} \partial_i(|w(\cdot,t)|^{r(t)-2})w_i\left(-\frac{\partial_i\partial_j}{\Delta}(v_iw_j+w_iv_j)\right) \mathrm{d}x,
$$
and
$$
\uppercase\expandafter{\romannumeral5}_2= r(t)^{2} \int_{\mathbb{R}^{3}} \partial_i(|w(\cdot,t)|^{r(t)-2})w_i\left(-\frac{\partial_i\partial_j}{\Delta}w_iw_j\right) \mathrm{d}x.
$$
According to \cite{Gra}, there holds $|x|^{r-2}\in A_r$ with $1<r<\infty$. By  H\"{o}lder's inequality, boundedness of the Riesz transforms on weighted $L^p$ spaces (Theorem 9.4.6 in \cite{Gra}), Lemma \ref{|x|v_c} and the Hardy inequality, there holds
\begin{eqnarray}\label{51}
\uppercase\expandafter{\romannumeral5}_1&\leq&2r(t)(r(t)-2)\int_{\mathbb{R}^{3}} \left| \nabla\left(|w(\cdot, t)|^{\frac{r(t)}{2}}\right)\right||w(\cdot, t)|^{\frac{r(t)}{2}-1} \left|\frac{\partial_i\partial_j}{\Delta}(v_iw_j+w_iv_j)\right|\mathrm{d} x\nonumber\\
&\leq& 4 r(t)(r(t)-2) C_{r}\||x|^{\frac{r-2}{r}}\left(v_{c} \otimes w\right)\|_{L^{r}}\left\| \frac{w^{\frac{r}{2}-1}}{x^{\frac{r-2}{r}}} \right\|_{L^{\frac{2 r}{r-2}}}\| \nabla(|w(\cdot, t)|^{\frac{r}{2}})\|_{L^{2}}\nonumber\\
&\leq& 4 r(t)(r(t)-2) C_{r}\||x|v_{c} \|_{L^{\infty}}\||x|^{-\frac{2}{r}} w\|_{L^{r}}\left\|\frac{w}{|x|^{\frac{2}{r}}}\right\|_{L^{r}}^{\frac{r}{2}-1}\| \nabla(|w(\cdot, t)|^{\frac{r}{2}})\|_{L^{2}}\nonumber\\
&\leq& 4 r(t)(r(t)-2) C_{r}K_c\left\|\frac{|w|^{\frac{r}{2}}}{|x|}\right\|_{L^{2}}\| \nabla(|w(\cdot, t)|^{\frac{r}{2}})\|_{L^{2}}\nonumber\\
&\leq& 8 r(t)(r(t)-2) C_{r} K_c\| \nabla(|w(\cdot, t)|^{\frac{r}{2}})\|_{L^{2}}^2,
\end{eqnarray}
where $C_r$ is as in Theorem 9.4.6 in \cite{Gra}.
Thanks to \cite{Iwa}, we deduce that $\|\frac{\partial_i\partial_j}{\Delta} f\|_{L^r} \leq H_r\|f\|_{L^r} $.
Combining with  H\"{o}lder's inequality and Sobolev embedding $\dot{H}^{1}(\mathbb{R}^{3})\hookrightarrow L^6(\mathbb{R}^{3})$, we have
\begin{eqnarray}
\uppercase\expandafter{\romannumeral5}_2&\leq & 2r(t)(r(t)-2) \|\nabla(|w(\cdot, t)|^{\frac{r(t)}{2}})\|_{L^2}\||w(\cdot, t)|^{\frac{r(t)}{2}-1} \|_{L^{\frac{6r}{r-2}}}\left\|\frac{\partial_i\partial_j}{\Delta}w_iw_j\right\|_{L^{\frac{3r}{r+1}}}\nonumber\\
&\leq& 2r(t)(r(t)-2) H_{\frac{3r}{r+1}}\|\nabla(|w(\cdot, t)|^{\frac{r(t)}{2}})\|_{L^2}\||w(\cdot, t)|^{\frac{r(t)}{2}-1} \|_{L^{\frac{6r}{r-2}}}\left\|w \otimes w\right\|_{L^{\frac{3r}{r+1}}}\nonumber\\
&\leq& 2r(t)(r(t)-2) H_{\frac{3r}{r+1}}\|\nabla(|w(\cdot, t)|^{\frac{r(t)}{2}})\|_{L^2}\|w(\cdot, t) \|^{\frac{r(t)}{2}-1}_{L^{3r}}\left\|w(\cdot, t)\right\|_{L^{3r}}\left\|w(\cdot, t)\right\|_{L^{3}}\nonumber\\
&\leq& 4r(t)(r(t)-2) H_{\frac{3r}{r+1}}\|\nabla(|w(\cdot, t)|^{\frac{r(t)}{2}})\|_{L^2}\||w(\cdot, t)|^{\frac{r(t)}{2}}\|_{L^6}\|w(\cdot, t) \|_{L^{3}}\nonumber\\
&\leq& 4r(t)(r(t)-2) H_{\frac{3r}{r+1}}\|\nabla(|w(\cdot, t)|^{\frac{r(t)}{2}})\|^2_{L^2}\|w(\cdot, t) \|_{L^{3}}.\nonumber
\end{eqnarray}{}
According to (\ref{L^3 stability-1}), when $\|w_0\|_{L^3}\leq \varepsilon_0,$ there holds
\begin{equation}\label{52}
\begin{aligned}
\uppercase\expandafter{\romannumeral5}_2 \leq 4r(t)(r(t)-2) H_{\frac{3r}{r+1}}C \varepsilon_0\|\nabla(|w(\cdot, t)|^{\frac{r(t)}{2}})\|^2_{L^2}.
\end{aligned}
\end{equation}
Combining with (\ref{important_1})-(\ref{52}), there holds
\begin{eqnarray}\label{remk_2}
&&r(t)^{2}\|w(t)\|_{r(t)}^{r(t)-1} \frac{\mathrm{d}}{\mathrm{d} t}\|w(t)\|_{r(t)}\nonumber\\
&\leq&\dot{r}(t) \int_{\mathbb{R}^{3}} |w(t)|^{r(t)} \ln \left(|w(t)|^{r(t)} /\|w(t)\|_{r(t)}^{r(t)}\right) \mathrm{d}x\nonumber\\
&&-\left(4r(t)(r(t)-2)\mu-2r(t)^2 K_c\right)\|\nabla (|w(t)|^{\frac{r(t)}{2}})\|_{L^2}^2,\nonumber
\end{eqnarray}
with
\begin{equation}\label{mu}
\mu=\inf_{3\leq r\leq q}\left\{1-\frac14C\varepsilon_0-\frac12 K_c-K_c-2C_rK_c-H_{\frac{3r}{r+1}}C \varepsilon_0\right\}> \frac 12.
\end{equation}
Applying the sharp logarithmic Sobolev inequality in \cite{Ca,Ngu}, we have
\begin{equation}
2\int|u|^2\ln\left(\frac{|u|}{\|u\|_{L^2}}\right)dx+3(1+\ln a)\|u\|_{L^2}^2\leq\frac{a^2}{\pi}\int |\nabla u|^2dx.
\end{equation}
Using $a=\left(\pi\frac{4r(t)(r(t)-2)\mu-2K_c r(t)^2}{ \dot{r}(t)}\right)^{\frac12}$ and $u=|w|^{\frac{r(t)}{2}}$, we obtain
\begin{eqnarray*}
&&r(t)^{2}\|w(t)\|_{r(t)}^{r(t)-1} \frac{\mathrm{d}}{\mathrm{d} t}\|w(t)\|_{r(t)}\\
&\leq& -\dot{r}(t)\left(3+\frac{3}{2} \ln \frac{(4\pi\mu -2\pi K_c) r(t)^2-8\pi \mu r(t) }{\dot{r}(t)}\right)\|w\|_{L^r}^r.
\end{eqnarray*}
If we define $G(t) :=\ln \left(\|w(t)\|_{L^{r(t)}}\right)$ with $r(t)=\frac{1}{\frac{1}{T}(\frac{1}{q}-\frac{1}{3})t+\frac{1}{3}}$, there holds
\begin{eqnarray*}
\frac{\mathrm{d} G(t)}{\mathrm{d} t} &\leq& -\frac{\dot{r}(t)}{r^{2}(t)}\left(3+\frac{3}{2} \ln\frac{(4 \pi\mu-2K_c \pi)r(t)^2-8\pi\mu r(t)}{\dot{r}(t)} \right)\nonumber\\
&=& \frac{1}{T}\left(\frac{1}{q}-\frac{1}{3}\right)\left(3+\frac{3}{2} \ln \left( \frac{-8\pi\mu}{r(t)}+4 \pi\mu-2K_c \pi \right) \right)
+\frac32\frac{1}{T}\left(\frac{1}{3}-\frac{1}{q}\right)\ln \frac{1}{T}\left(\frac{1}{3}-\frac{1}{q}\right)\nonumber\\
&\leq& \frac32\frac{1}{T}\left(\frac{1}{3}-\frac{1}{q}\right)\ln \frac{1}{T}\left(\frac{1}{3}-\frac{1}{q}\right).\nonumber
\end{eqnarray*}
Integrating this in time from 0 to $T$ yields
$$
G(T) \leq  G(0)+ \frac32\left(\frac1q-\frac13\right)\ln \left(\frac{T}{1 / 3-1 / q}\right).
$$
We obtain
$$
\ln \|w(\cdot, T)\|_{r(T)}\leq \ln \|w(\cdot, 0)\|_{L^3}+\frac32\left(\frac1q-\frac13\right)\ln \left(\frac{T}{1 / 3-1 / q}\right).
$$
Hence, we obtain
\begin{equation}\label{w}
\|w(\cdot, t)\|_{L^q}\leq  C_q t^{\frac32\left(\frac1q-\frac13\right)}\|w(\cdot, 0)\|_{L^3},
\end{equation}
with $C_q=(\frac13-\frac1q)^{\frac32 (\frac13-\frac1q)}.$ 

To give strict proof of (\ref{nonlinear decay}),  we consider the  approximation scheme. Using method in \cite{Kw,Lem1,Lem},  the mollified system in $\mathbb{R}^{3} \times(0, \infty)$ is as follows
\begin{equation}\label{PNS^e}
\begin{cases}
 w^\epsilon_{t}-\Delta w^\epsilon+\left(\mathcal{J}_{\epsilon}\left(w^{\epsilon}\right) \cdot \nabla\right)w^{\epsilon}+(\mathcal{J}_{\epsilon}(w^\epsilon) \cdot \nabla) v_{c}+\left(v_{c} \cdot \nabla\right)\mathcal{J}_{\epsilon}(w^\epsilon) +\nabla \pi^\epsilon=0,\\
 \nabla\cdot w^\epsilon=0,\\
 w^\epsilon(x, 0) =w_{0}(x),
\end{cases}
\end{equation}
where $\mathcal{J}_{\epsilon}(v)=v * \eta_{\epsilon}, \epsilon>0,$ the mollifier $\eta_{\epsilon}(x)=\varepsilon^{-3} \eta\left(\frac{x}{\epsilon}\right)$ with positive $\eta\in C_c^{\infty}(B(0,1)),$ $\int \eta dx=1.$
By classical approximation method,  solution $w^\epsilon$ satisfies (\ref{nonlinear decay}).
Similar as (\ref{remk_1}), we have
\begin{eqnarray*}
&&r(t)^{2}\|w^\epsilon(\cdot,t)\|_{r(t)}^{r(t)-1} \frac{\mathrm{d}}{\mathrm{d} t}\|w^\epsilon(\cdot,t)\|_{r(t)}\\
&=&\dot{r}(t) \int_{\mathbb{R}^{3}} |w^\epsilon(\cdot,t)|^{r(t)} \ln \left(|w^\epsilon(\cdot,t)|^{r(t)} /\|w^\epsilon(\cdot,t)\|_{r(t)}^{r(t)}\right) \mathrm{d} x\\
&&+r(t)^{2} \int_{\mathbb{R}^{3}} |w^\epsilon(\cdot,t)|^{r(t)-2} w^\epsilon_i\Big( \partial_j\partial_j w^\epsilon_i-\partial_j (\mathcal{J}_{\epsilon}(w^\epsilon)_iw^\epsilon_j
  +\mathcal{J}_{\epsilon}(w^\epsilon)_iv_j+v_i\mathcal{J}_{\epsilon}(w^\epsilon)_j)-\partial_i\pi^\epsilon\Big)  \mathrm{d}x\\
&:=&\dot{r}(t) \int_{\mathbb{R}^{3}} |w^\epsilon(\cdot,t)|^{r(t)} \ln \left(|w^\epsilon(\cdot,t)|^{r(t)} /\|w^\epsilon(\cdot,t)\|_{r(t)}^{r(t)}\right) \mathrm{d} x\\
&& +\uppercase\expandafter{\romannumeral1}^\epsilon+\uppercase\expandafter{\romannumeral2}^\epsilon +\uppercase\expandafter{\romannumeral3}^\epsilon+\uppercase\expandafter{\romannumeral4}^\epsilon+\uppercase\expandafter{\romannumeral5}^\epsilon.
\end{eqnarray*}
Integration by parts show that
\begin{equation}
\begin{aligned}
\uppercase\expandafter{\romannumeral1}^\epsilon&= -r(t)^{2} \int_{\mathbb{R}^{3}} \partial_j(|w^\epsilon(\cdot, t)|^{r(t)-2} w^\epsilon_i) \partial_j w^\epsilon_idx \\
&= -4r(t)(r(t)-2)\|\nabla (|w^\epsilon(\cdot,t)|^{\frac{r(t)}{2}})\|_{L^2}^2-r(t)^{2} \int_{\mathbb{R}^{3}}  |w^\epsilon(t)|^{r(t)-2}|\nabla w^\epsilon|^2dx.
\end{aligned}
\end{equation}
Similar as (\ref{important_1})-(\ref{52}), we have
\begin{eqnarray}
\uppercase\expandafter{\romannumeral2}^\epsilon&=& -{r(t)^{2} }\int_{\mathbb{R}^{3}} |w^\epsilon(\cdot,t)|^{r(t)-2} w^\epsilon_i\partial_j (\mathcal{J}_{\epsilon}(w^\epsilon)_iw^\epsilon_j\mathrm{d}x\nonumber\\
&=& {r(t)^{2} }\int_{\mathbb{R}^{3}}\partial_j(|w^\epsilon(\cdot,t)|^{r(t)-2}) w^\epsilon_i\mathcal{J}_{\epsilon}(w^\epsilon)_iw^\epsilon_j\mathrm{d}x
  +{r(t)^{2} }\int_{\mathbb{R}^{3}} |w^\epsilon(\cdot,t)|^{r(t)-2}\partial_j( w^\epsilon_i)\mathcal{J}_{\epsilon}(w^\epsilon)_iw^\epsilon_j\mathrm{d}x\nonumber\\
&\leq&  r(t)(r(t)-2)C\varepsilon_0 \left\|\nabla\left(|w^\epsilon(t)|^{\frac{r(t)}{2}}\right)\right\|_{L^{2}}^{2},
\end{eqnarray}
\begin{equation}
\uppercase\expandafter{\romannumeral3}^\epsilon\leq 2 r(t)(r(t)-2) K_{c}\left\|\nabla\left(|w^\epsilon(t)|^{\frac{r(t)}{2}}\right)\right\|_{L^{2}}^{2},
\end{equation}
\begin{eqnarray}
\uppercase\expandafter{\romannumeral4}^\epsilon&\leq&  \frac{r(t)^{2}}{2} K_c\int_{\mathbb{R}^{3}} |w^\epsilon(\cdot,t)|^{r(t)-2} |\nabla w^\epsilon(\cdot,t)|^{2} dx\\
&&+ (2 r(t)(r(t)-2)+4r(t)^{2}) K_c\left\|\nabla\left(|w^\epsilon(t)|^{\frac{r(t)}{2}}\right)\right\|_{L^{2}}^{2},\nonumber
\end{eqnarray}
\begin{equation}
\uppercase\expandafter{\romannumeral5}^\epsilon\leq 4 r(t)(r(t)-2)(2C_{r} K_c+H_{\frac{3r}{r+1}}C \varepsilon_0) \| \nabla(|w^\epsilon(\cdot, t)|^{\frac{r}{2}})\|_{L^{2}}^2,
\end{equation}
and
\begin{eqnarray}
&&r(t)^{2}\|w^\epsilon(t)\|_{r(t)}^{r(t)-1} \frac{\mathrm{d}}{\mathrm{d} t}\|w^\epsilon(t)\|_{r(t)}\\
&\leq& \dot{r}(t) \int_{\mathbb{R}^{3}} |w^\epsilon(t)|^{r(t)} \ln \left(|w^\epsilon(t)|^{r(t)} /\|w^\epsilon(t)\|_{r(t)}^{r(t)}\right) \mathrm{d}x\nonumber\\
&&-\left(4r(t)(r(t)-2)\mu-2r(t)^2 K_c\right)\|\nabla (|w^\epsilon(t)|^{\frac{r(t)}{2}})\|_{L^2}^2,\nonumber
\end{eqnarray}
with
\begin{equation}
\mu=\inf_{3\leq r\leq q}\left\{1-\frac14C\varepsilon_0-\frac12 K_c-K_c-2C_rK_c-H_{\frac{3r}{r+1}}C \varepsilon_0\right\}> \frac 12.
\end{equation}
Similar as procedure  in the proof of (\ref{w}), we obtain
\begin{equation}
\|w^\epsilon(\cdot, t)\|_{L^q}\leq  (\frac13-\frac1q)^{\frac32 (\frac13-\frac1q)}t^{\frac32\left(\frac1q-\frac13\right)}\|w^\epsilon(\cdot, 0)\|_{L^3}.
\end{equation}
By compactness and convergence theory,  solution $w$ satisfies (\ref{nonlinear decay}).

Finally, we will prove (\ref{L^3 stability-2}). Since $w_0\in L^3$ and $\|w_0\|_{L^3}<\varepsilon_0$, there exists a subsequence denoted by $\{w_{0,n}\}$ such that $w_{0,n}\in L^2\cap L^3$  and
$$
w_{0,n}\rightarrow w_0 \text{\ in\ }L^3 \text{\ as } n\rightarrow \infty.
$$
According to Theorem \ref{s-w-uniqueness}, Corollary \ref{corollary} and Remark \ref{remark-Kar}, we have
$$
\lim_{t\rightarrow \infty}\|w_n(\cdot, t)\|_{L^3}=0.
$$
Based on similar  proof of (\ref{Z-crucial}), when $\|w_{0,n}-w_0\|_{L^3}\leq (4C^2e^{2C \int_0^T \|w\|_{L^6}^4 dt})^{-1},$ we have
$$
\|w_n-w\|_{L_t^{\infty}([0,\infty);L_x^3)}\leq 2C\|w_{0,n}-w_0\|_{L^3}e^{C\int_0^{\infty} \|w\|_{L^6}^4 dt},
$$
for a  positive constant $C.$ According to Theorem \ref{L^3 well-posedness},  we have
$$
\int_0^{\infty} \|w\|_{L^6}^4 dt\leq C \|w_0\|_{L^3(\mathbb{R}^3)}.
$$
Hence,
$$
\lim_{n\rightarrow \infty}\|w_n-w\|_{L_t^{\infty}([0,\infty);L_x^3)}=0.
$$
Therefore, (\ref{L^3 stability-2}) holds.
 {\hfill
$\square$\medskip}

\section{The linear operator $\mathcal{L}$ }\label{linear}

By Lemmas \ref{alemm} and \ref{p_alemm}, as stated in Sections  \ref{Sec2} and \ref{Sec4},
$-\mathcal{L}$ is the  infinitesimal generator of the strongly continuous  semigroup of contraction $e^{-t\mathcal{L}}$ of  bounded linear operators on  $L_{\sigma}^q(\mathbb{R}^3)$, $1< q<\infty$. In this section, we prove that $e^{-t\mathcal{L}}$  is an analytic semigroup.

We consider the following system
 \begin{equation}\label{lambda-PNS}
\begin{cases}
 \lambda u -\Delta u+(u \cdot \nabla) v_{c}+\left(v_{c} \cdot \nabla\right) u+\nabla p=f,\\
 \nabla\cdot u=0.
\end{cases}
\end{equation}
For $\delta>0$ small, set $\Sigma_{\delta}=\{  \lambda\in \mathbb{C}\backslash \{ 0 \}: |\arg \lambda|<\frac{\pi}{2}+ \delta\}$. It is easy to see that for $\lambda = \sigma +\sqrt{-1}\tau  \in \Sigma,$ $\sigma$, $\tau$ real, if $\sigma <0,$ then
\begin{equation}\label{2-11}
|\sigma |<\delta|\tau|.
\end{equation}

\begin{theorem}\label{linear ope property}
For $1<q<\infty,$ there exist some positive  constants $\delta$ and $\bar{c}_q$ which depend only on $q$  such that for any $|c|>\bar{c}_q,$ $\lambda \in \Sigma_{\delta},$ and $u\in C_{c,\sigma}^{\infty}(\mathbb{R}^3)$ satisfying system (\ref{lambda-PNS}), we have
\begin{equation}\label{2-1}
\|u\|_{L^q}\leq \frac{C}{|\lambda|}\|f\|_{L^q},
\end{equation}
where $C$ is a constant depending only on $q$ and $\delta.$
Consequently, $e^{-t\mathcal{L}}$ is an  analytic semigroup of bounded linear operators on $L^q_{\sigma}(\mathbb{R}^3)$ in the sector $\{\lambda: |arg \lambda|<\delta\}$.
\end{theorem}
The last statement in the above theorem follows from estimate (\ref{2-1}), together with the fact that $e^{-t\mathcal{L}}$ is a
 strongly continuous  semigroup  of contraction on  $L_{\sigma}^q(\mathbb{R}^3)$ for  $1< q<\infty$ which are established in Sections \ref{Sec2} and  \ref{Sec4}, see Theorem 1.5.2 in \cite{Paz}.

The following theorem for operators $\mathcal{L} u$ on scalar functions $u$ can be proved by using the arguments in the proof of Theorem \ref{linear ope property}, but much simpler since no pressure term is present.

\begin{theorem}
For $n\geq 3$ and $1<q<\infty,$ there exist some positive constants $\varepsilon$ and $\delta$ which depend only on $n$ and $q$ such that the operator
$$\mathcal{L}u:=-\Delta u+a(x)u+b(x)\cdot  \nabla u,$$
 with  $|a(x)|\leq \varepsilon|x|^{-2}$ and $|b(x)|\leq \varepsilon|x|^{-1}$ for all $x\in \mathbb{R}^n,$ has the property that
 $e^{-t\mathcal{L}}$ is an  analytic semigroup of bounded linear operators on $L^q_{\sigma}(\mathbb{R}^3)$ in the sector $\{\lambda: |arg \lambda|<\delta\}$.
\end{theorem}

To prove Theorem \ref{linear ope property}, we need the following lemma.
\begin{lemma}\label{sec3-lem1}
Function $u$ has the following property
\begin{equation}
|\nabla(|u|^2)|^2\leq 4|\nabla u|^2|u|^2.
\end{equation}
Consequently, for $1\leq q\leq2,$ we have
\begin{equation}
\frac{q-2}{4}\int_{\mathbb{R}^3}|u|^{q-4} |\nabla (|u|^2)|^2dx+ \int_{\mathbb{R}^3} |\nabla u|^2|u|^{q-2}dx \geq (q-1)\int_{\mathbb{R}^3} |\nabla u|^2|u|^{q-2}dx.
\end{equation}
\end{lemma}
\noindent\textbf{Proof.}
We have
\begin{equation}
|\nabla(|u|^2)|^2=\sum_j|\partial_j(|u|^2)|^2=\sum_j|\partial_j<u,u>|^2\leq4\sum_j|<\partial_ju,u>|^2.
\end{equation}
For fixed $j,$ using Cauchy-Schwartz inequality, we have
\begin{eqnarray}
|\nabla(|u|^2)|^2\leq 4\sum_j(|\partial_j u|^2|u|^2)=4|\nabla u|^2|u|^2.
\end{eqnarray}
 {\hfill
$\square$\medskip}

\noindent\textbf{Proof of Theorem \ref{linear ope property}.}
The value of $\delta$ will be chosen in the proof below.
Multiplying the equation ($\ref{lambda-PNS})_1$ by $|u|^{q-2}\overline{u}$, and integrating it on $\mathbb{R}^3$, we have
\begin{equation}\label{Sec3_0}
\begin{aligned}
\int_{\mathbb{R}^3} \nabla u \cdot \nabla(|u|^{q-2}\overline{u}) dx+ \lambda \int_{\mathbb{R}^3} |u|^q dx+ \int_{\mathbb{R}^3} (v_c \cdot \nabla u )\cdot (|u|^{q-2}\overline{u}) dx&\\
+ \int_{\mathbb{R}^3} (u \cdot \nabla v_c )\cdot (|u|^{q-2}\overline{u}) dx+\int_{\mathbb{R}^3} \nabla p \cdot (|u|^{q-2}\overline{u}) dx
&= \int_{\mathbb{R}^3} f \cdot (|u|^{q-2}\overline{u}) dx.
\end{aligned}
\end{equation}
Set $\uppercase\expandafter{\romannumeral 1}+\uppercase\expandafter{\romannumeral 2}+\uppercase\expandafter{\romannumeral 3}+\uppercase\expandafter{\romannumeral 4}+\uppercase\expandafter{\romannumeral 5} =\int_{\mathbb{R}^3} f \cdot (|u|^{q-2}\overline{u} )dx. $
For the first part, we have
\begin{eqnarray*}
\uppercase\expandafter{\romannumeral 1} &=&\int_{\mathbb{R}^3} \partial_j u_i \partial_j\left((u_m \overline{u}_m)^{\frac{q-2}{2}} \overline{u}_i\right)  dx \\
&=&\int_{\mathbb{R}^3} (\partial_j u_i) \overline{u}_i \frac{q-2}{2}|u|^{q-4}\partial_j(|u|^2)  dx + \int_{\mathbb{R}^3}  (\partial_j u_i)(\overline{\partial_j u_i})|u|^{q-2} dx\\
&=&\int_{\mathbb{R}^3} (\partial_j u_i) \overline{u}_i \frac{q-2}{2}|u|^{q-4}\left[(\partial_ju_m) \overline{u}_m+ \overline {(\partial_ju_m) \overline{u}_m}\right]  dx\\
&& + \int_{\mathbb{R}^3}  (\partial_j u_i)(\overline{\partial_j u_i})|u|^{q-2} dx\\
&=&\int_{\mathbb{R}^3} (\partial_j u_i) \overline{u}_i \frac{q-2}{2}|u|^{q-4}\left[(\partial_ju_m) \overline{u}_m+ \overline {(\partial_ju_m) \overline{u}_m}\right]  dx + \int_{\mathbb{R}^3}  |\nabla u|^2|u|^{q-2} dx.
\end{eqnarray*}
Denote $\xi_j =(\partial_j u_m) \overline{u}_m=a_j+\sqrt{-1}b_j, $ then
\begin{eqnarray}
\partial_j(|u|^2)= (\partial_j u_m) \overline{u}_m+ \overline {(\partial_j u_m) \overline{u}_m}=2 \text{ Re }\xi_j= 2 a_j.
\end{eqnarray}
Therefore, we have
\begin{eqnarray*}
\uppercase\expandafter{\romannumeral 1}
&=&\int_{\mathbb{R}^3}  \frac{q-2}{2}|u|^{q-4}\xi_j\left[\xi_j+ \overline{\xi}_j\right]   +  |\nabla u|^2|u|^{q-2} dx\\
&=&\int_{\mathbb{R}^3}  \frac{q-2}{2}|u|^{q-4}(a_j+\sqrt{-1}b_j)\left(2a_j\right)   +  |\nabla u|^2|u|^{q-2} dx\\
&=&\int_{\mathbb{R}^3}(q-2)|u|^{q-4} \Sigma(a_j)^2 dx+  \sqrt{-1} \int_{\mathbb{R}^3}\frac{q-2}{2}|u|^{q-4}  2a_j b_j dx
 + \int_{\mathbb{R}^3} |\nabla u|^2|u|^{q-2} dx\\
&=&\int_{\mathbb{R}^3}\frac{q-2}{4}|u|^{q-4} |\nabla (|u|^2)|^2 dx+  \sqrt{-1} \int_{\mathbb{R}^3}\frac{q-2}{2}|u|^{q-4}  2a_j b_j dx
+ \int_{\mathbb{R}^3} |\nabla u|^2|u|^{q-2} dx.
\end{eqnarray*}
For the second part $\uppercase\expandafter{\romannumeral 2}$
\begin{eqnarray}
\uppercase\expandafter{\romannumeral 2} = (\sigma +\sqrt{-1}\tau)  \int_{\mathbb{R}^3} |u|^q dx.
\end{eqnarray}
It is easy to see that
\begin{eqnarray}\label{ineqaulity}
2 \sum_j|a_j||b_j|\leq |\nabla u|^2|u|^2.
\end{eqnarray}
Then, using Lemma \ref{sec3-lem1} and (\ref{ineqaulity}), we have
\begin{eqnarray}\label{ineq3_1}
|\uppercase\expandafter{\romannumeral 1}+\uppercase\expandafter{\romannumeral 2} |&\geq &
\text{Re }(\uppercase\expandafter{\romannumeral 1}+\uppercase\expandafter{\romannumeral 2}) +|\text{Im }(\uppercase\expandafter{\romannumeral 1}+\uppercase\expandafter{\romannumeral 2}) |\nonumber\\
&\geq&
 \text{min}\{q-1,1\}\int_{\mathbb{R}^3} |\nabla u|^2|u|^{q-2}dx+\sigma\int_{\mathbb{R}^3}|u|^{q} dx\nonumber
 +\left|\tau\int_{\mathbb{R}^3}|u|^{q} dx-\int_{\mathbb{R}^3} \frac{q-2}{2}|u|^{q-4}2a_jb_jdx\right|\nonumber\\
&\geq&
\text{min}\{q-1,1\}\int_{\mathbb{R}^3} |\nabla u|^2|u|^{q-2}dx+\sigma\int_{\mathbb{R}^3}|u|^{q} dx\nonumber\\
&&+\left||\tau|\int_{\mathbb{R}^3}|u|^{q} dx-\left|\frac{q-2}{2}\right|\int_{\mathbb{R}^3} |u|^{q-4}|\nabla u|^2|u|^{q-2}dx\right|.
\end{eqnarray}
We distinguish into two cases:

\noindent\textbf{Case 1.} $ \text{min}\{q-1,1\}\int_{\mathbb{R}^3} |\nabla u|^2|u|^{q-2} dx\geq 8\delta|\tau| \int_{\mathbb{R}^3}|u|^{q} dx.$

\noindent\textbf{Case 2.} $ \text{min}\{q-1,1\}\int_{\mathbb{R}^3} |\nabla u|^2|u|^{q-2} dx< 8\delta|\tau| \int_{\mathbb{R}^3}|u|^{q} dx.$

In Case 1, we deduce from (\ref{ineq3_1}), using (\ref{2-11}) and requiring $0<\delta\leq 1,$ that
\begin{eqnarray}
|\uppercase\expandafter{\romannumeral 1}+\uppercase\expandafter{\romannumeral 2} |&\geq\frac12\text{min}\{q-1,1\}\int_{\mathbb{R}^3} |\nabla u|^2|u|^{q-2} dx+(4\delta|\tau|+\sigma) \int_{\mathbb{R}^3}|u|^{q} dx\nonumber\\
&\geq \frac12\text{min}\{q-1,1\}\int_{\mathbb{R}^3} |\nabla u|^2|u|^{q-2} dx+\delta(|\tau|+|\sigma|) \int_{\mathbb{R}^3}|u|^{q} dx.
\end{eqnarray}

In Case 2, we derive from  (\ref{ineq3_1}) that
\begin{equation}
|\uppercase\expandafter{\romannumeral 1}+\uppercase\expandafter{\romannumeral 2} |\geq \text{min}\{q-1,1\}\int_{\mathbb{R}^3} |\nabla u|^2|u|^{q-2} dx+(|\tau|+\sigma-\frac{4|q-2|\delta |\tau|}{\text{min}\{q-1,1\}}) \int_{\mathbb{R}^3}|u|^{q} dx.
\end{equation}
Now we require that $\delta$ further satisfies $8|q-2|\delta<\text{min}\{q-1,1\}$ and $\delta\leq \frac 15.$ Then we have
\begin{eqnarray*}
|\uppercase\expandafter{\romannumeral 1}+\uppercase\expandafter{\romannumeral 2} |&\geq& \text{min}\{q-1,1\}\int_{\mathbb{R}^3} |\nabla u|^2|u|^{q-2} dx+(\frac 12|\tau|+\sigma) \int_{\mathbb{R}^3}|u|^{q} dx\\
&\geq& \text{min}\{q-1,1\}\int_{\mathbb{R}^3} |\nabla u|^2|u|^{q-2} dx+\frac 14(|\tau|+\sigma) \int_{\mathbb{R}^3}|u|^{q} dx.
\end{eqnarray*}
So in both cases, we have proved that
\begin{equation}\label{Sec3_12}
|\uppercase\expandafter{\romannumeral 1}+\uppercase\expandafter{\romannumeral 2} |\geq \frac12\text{min}\{q-1,1\}\int_{\mathbb{R}^3} |\nabla u|^2|u|^{q-2} dx+\delta(|\tau|+|\sigma|) \int_{\mathbb{R}^3}|u|^{q} dx.
\end{equation}
By  H\"{o}lder's inequality, Hardy inequality, Cauchy inequality and Lemma \ref{sec3-lem1},  we have
\begin{eqnarray}\label{Sec3_3}
|\uppercase\expandafter{\romannumeral 3}|&\leq& K_c \int_{\mathbb{R}^3} \frac{|\nabla u|}{|x|} |u|^{q-1}dx\nonumber\\
&\leq& K_c  \left\| |\nabla u||u|^{\frac{q-2}{2}}\right\|_{L^2}\left\|\frac{|u|^{\frac q2}}{|x|}\right\|_{L^2}\nonumber\\
&\leq& 2K_c  \left\| |\nabla u||u|^{\frac{q-2}{2}}\right\|_{L^2}\left\|\nabla(|u|^{\frac q2})\right\|_{L^2}\\
&\leq& 2K_c \int_{\mathbb{R}^3} |\nabla u|^2 |u|^{q-2} dx+ CK_c \int_{\mathbb{R}^3} |u|^{q-4}|\nabla(| u|^2)|^2 dx\nonumber\\
&\leq&  CK_c \int_{\mathbb{R}^3} |\nabla u|^2 |u|^{q-2} dx.\nonumber
\end{eqnarray}
Similarly,  by Hardy inequality and  Lemma \ref{sec3-lem1},  we deduce
\begin{eqnarray}\label{Sec3_4}
|\uppercase\expandafter{\romannumeral 4}|&\leq& K_c \int_{\mathbb{R}^3} \frac{|u|^{q}}{|x|^2}dx\nonumber\\
&\leq& 2 K_c \int_{\mathbb{R}^3} \nabla |(|u|^{\frac q2})|^2dx\\
&\leq&  CK_c \int_{\mathbb{R}^3} |\nabla u|^2 |u|^{q-2} dx.\nonumber
\end{eqnarray}
According to (\ref{lambda-PNS}), we have
\begin{eqnarray*}
p= \frac{\text{div}}{\Delta}f-\frac{\partial_i\partial_j}{\Delta}(v_c \otimes u+ u \otimes v_c)_{ij}.
\end{eqnarray*}
 By integration by parts, we have
\begin{equation}\label{p}
\uppercase\expandafter{\romannumeral 5}=\int_{\mathbb{R}^3} \frac{\partial_i\partial_j}{\Delta}f \cdot(|u|^{q-2} \overline{u}) dx +\int_{\mathbb{R}^3} \frac{\partial_i\partial_j}{\Delta}(v_c \otimes u+ u \otimes v_c)_{ij} \nabla\cdot(|u|^{q-2} \overline{u}) dx.
\end{equation}
By boundedness of the Riesz transforms on weighted $L^p$ spaces (Theorem 9.4.6 in \cite{Gra}) and  H\"{o}lder's inequality, we have the following estimate
\begin{eqnarray}
\int_{\mathbb{R}^3} \frac{\partial_i\partial_j}{\Delta}f \cdot(|u|^{q-2} \overline{u}) dx\leq C\| f\|_{L^q}\| u\|_{L^q}^{q-1}.
\end{eqnarray}
By  H\"{o}lder's inequality, Hardy inequality, Sobolev embedding, boundedness of the Riesz transforms on weighted $L^p$ spaces (Theorem 9.4.6 in \cite{Gra}) and  Lemma \ref{sec3-lem1}, we have
\begin{eqnarray}\label{p_2}
&&\int_{\mathbb{R}^3} \frac{\partial_i\partial_j}{\Delta}(v_c \otimes u+ u \otimes v_c)_{ij} \nabla\cdot(|u|^{q-2} \overline{u}) dx \nonumber\\
&\leq& C\||x|^{\frac{q-2}{q}}\left(u\otimes v_c+v_c\otimes u\right)\|_{L^{q}} \|\nabla (|u|^{\frac q2})\|_{L^{2}} \left\|\frac{|u|^{\frac q2-1}}{|x|^{\frac{q-2}{q}}}\right\|_{L^{\frac{2q}{q-2}}}\nonumber\\
&\leq& C\||x|v_c\|_{L^{\infty}} \left\|\frac{u}{|x|^{\frac2q}}\right\|_{L^q}\|\nabla (|u|^{\frac q2})\|_{L^{2}} \left\|\nabla(|u|^{\frac q2})\right\|_{L^{2}}^{\frac {q-2} {q} }\\
&\leq& C K_c\left\|\frac{u}{|x|^{\frac2q}}\right\|_{L^q}\|\nabla (|u|^{\frac{q}{2}})\|_{L^2}^{\frac{2q-2}{q}}\nonumber\\
&\leq& C K_c \|\nabla (|u|^{\frac{q}{2}})\|^2_{L^2}\nonumber\\
&\leq& CK_c \int_{\mathbb{R}^3} |\nabla u|^2 |u|^{q-2} dx.\nonumber
\end{eqnarray}
Combining with (\ref{p})-(\ref{p_2}), we deduce
\begin{eqnarray}\label{Sec3_5}
|\uppercase\expandafter{\romannumeral 5}|\leq  C\| f\|_{L^q}\| u\|_{L^q}^{q-1}+CK_c \int_{\mathbb{R}^3} |\nabla u|^2 |u|^{q-2} dx.
\end{eqnarray}
Combining with (\ref{Sec3_0}), (\ref{Sec3_12})-(\ref{Sec3_4}), (\ref{Sec3_5}) and the condition of $K_c$ small enough, by  H\"{o}lder's inequality, we have
\begin{equation}
\begin{aligned}
 \frac12\text{min}\{q-1,1\}\int_{\mathbb{R}^3} |\nabla u|^2|u|^{q-2} dx+\delta(|\tau|+|\sigma|) \int_{\mathbb{R}^3}|u|^{q} dx\leq  C\| f\|_{L^q}\| u\|_{L^q}^{q-1}.
\end{aligned}
\end{equation}
Since $\lambda=\sigma+\sqrt{-1}\tau$, we deduce
\begin{equation}
\begin{aligned}
C(\delta)|\lambda|\| u\|_{L^q}\leq C\| f\|_{L^q},
\end{aligned}
\end{equation}
which implies (\ref{2-1}).
 {\hfill
$\square$\medskip}


\section{Weak-Strong Uniqueness}\label{w-s}
In this section, we will prove the Theorem \ref{s-w-uniqueness}, Proposition \ref{prop-s-w-uniqueness} and illustrate Corollary \ref{corollary} briefly.

First, we give the detailed proof of Theorem \ref{s-w-uniqueness}.

\noindent\textbf{Proof of Theorem \ref{s-w-uniqueness}.}
Following the proof of Theorem 4.4 in \cite{Tsa},  setting $g=v-u$, we have
\begin{equation}\label{f}
\begin{cases}
 \partial_t g-\Delta g+\nabla \pi=-((u+g) \cdot \nabla) g-(g\cdot \nabla) u-g\cdot \nabla v_c-v_c\cdot \nabla g,\\
 \nabla\cdot g=0,\\
 g(x, 0) =0.
\end{cases}
\end{equation}
Using $g$ itself as a test function and  integrating this in time from 0 to $t$, we have
\begin{equation}
\int \frac{|g|^{2}}{2} d x+\int_0^t\int|\nabla g|^{2} d xdt\leq\int_0^t\int (u+v_c) \cdot(g \cdot \nabla) g d xdt.
\end{equation}
Denote $E(t)=\operatorname{ess} \sup _{s<t}\|g(s)\|_{2}^{2}+\int_{0}^{t}\|\nabla g\|_{2}^{2}d \tau$ and $t_{0}=\sup \{t \in[0,T] : g(s)=0 \text { if } 0<s<t\}$.

We claim that $t_0=T.$ Using contradiction argument, we assume that $t_0<T.$
Since
\begin{equation}
\left|\iint u v \nabla w d x d t \right| \leq C\|u\|_{L_{t}^{s} L_{x}^{q}}\|v\|_{L_{t}^{\infty} L_{x}^{2}}^{2 / s}\|v\|_{L_{t}^{2} L_{x}^{6}}^{3 / q}\|\nabla w\|_{L_{t, x}^{2}},
\end{equation}
for $\frac 3q+ \frac 2s=1$ with $1\leq q,s \leq \infty,$ we have
\begin{equation}
\left|\int_{t_{0}}^{t} \int u \cdot(g \cdot \nabla) g d x d \tau\right|\leq C\|u\|_{L_{t}^{s} L_{x}^{q}}E(t),
\end{equation}
for $t\in [t_0,T,]$
By H\"{o}lder inequality, Hardy inequality and Lemma \ref{|x|v_c}, we have
\begin{eqnarray}
\left|\int_{t_{0}}^{t} \int v_c \cdot(g \cdot \nabla) g d x d \tau\right|
&\leq& \int_{t_{0}}^{t} \||x|v_c\|_{L^\infty}\left\|\frac{g}{|x|}\right\|_{L^2}\|\nabla g\|_{L^2}d \tau\nonumber\\
&\leq& 2\int_{t_{0}}^{t} \||x|v_c\|_{L^\infty}\|\nabla g\|^2_{L^2}d \tau\nonumber\\
&\leq& 2K_c\|\nabla g\|_{L_{t,x}^2}^2\nonumber\\
&\leq& 2K_cE(t),
\end{eqnarray}
for $t\in [t_0,T].$
Hence, there holds
\begin{equation}
\begin{aligned}
E(t)\leq C\|u\|_{L_{t}^{s}([t_0,t]; L_{x}^{q})}E(t)+2 K_cE(t).
\end{aligned}
\end{equation}
If $s<\infty,$ we have $C\|u\|_{L_{t}^{s}([t_0,t]; L_{x}^{q})}<\frac 14$ for $t$ sufficiently close to $t_0$. If $s=\infty,$ we need $C\|u\|_{L_{t}^{\infty}([t_0,t]; L_{x}^{3})}<\frac 14$. Moreover, $K_c<\frac14$ by assumption. Therefore, $E(s)=0$  for all $s\in[t_0,t],$ which makes a contradiction to the definition of $t_0$. Hence, $t_0=T$ and  for all $t\in [0,T]$.
 {\hfill
$\square$\medskip}

Following is the proof of Proposition \ref{prop-s-w-uniqueness}.

\noindent\textbf{Proof of Proposition \ref{prop-s-w-uniqueness}.}
For $p\geq 3,$ our goal is to show that the $L^p$ mild solution $w$ is a $L^2$-weak solution. Crucial part is to prove that $w\in   C_w([0,T];L_{x}^2)\cap L_T^2 (\dot{H}_{x}^1)$. Set $w=a+z$ as in Section \ref{Sec2}.  We will prove $a\in   C_w([0,T];L_{x}^2)\cap L_T^2 (\dot{H}_{x}^1)$ and $z\in  C_w([0,T];L_{x}^2)\cap L_T^2 (\dot{H}_{x}^1)$ as follows.

Multiplying $(\ref{a})_1$ by $a,$ then integrating it on $\mathbb{R}^3$, we have
\begin{equation}\label{4.7}
\frac{1}{2}\frac{\mathrm{d}}{\mathrm{d} t}\|a(t)\|^2_{L^2}+\|\nabla a\|^2_{L^2}
=-\int_{\mathbb{R}^3}\text{div}(a\otimes v_c+v_c\otimes a)\cdot a dx-\int_{\mathbb{R}^3} \nabla \pi_1 \cdot a dx.
\end{equation}
By similar estimate as (\ref{prop- rh1}), using integration by parts and div$a=0$, we obtain
\begin{equation}
-\int_{\mathbb{R}^3}\text{div}(a\otimes v_c+v_c\otimes a)\cdot a dx \leq 2K_c\left\|\nabla a\right\|_{L^2}^2,
\end{equation}
and
\begin{equation}\label{4.9}
-\int_{\mathbb{R}^3} \nabla \pi_1 \cdot a dx=0.
\end{equation}
 From (\ref{4.7})-(\ref{4.9}), we have
\begin{equation}
\frac 12\frac{\mathrm{d}}{\mathrm{d} t}\|a(t)\|^2_{L^2}+\|\nabla a\|^2_{L^2}\leq 2K_c  \|\nabla a\|^2_{L^2}.
\end{equation}
Since $|c|>{c}_p$ where ${c}_p$ is as in Theorem \ref{p>3 result}, we can guarantee $1-2 K_c>0.$ Combining with similar argument as (\ref{con_1}), we have $a\in  C([0,T];L_{x}^2)\cap L_T^2 (\dot{H}_{x}^1)$.

When $p\in[3,4]$, using $(w_1,w_2)=(w,w),$ multiplying $(\ref{N})_1$ by $z$ and integrating it on $\mathbb{R}^3$, we obtain
\begin{eqnarray}
&&\frac{1}{2}\frac{\mathrm{d}}{\mathrm{d} t}\|z(t)\|^2_{L^2}+\|\nabla z\|^2_{L^2}\nonumber\\
&=&-\int_{\mathbb{R}^3}\text{div}(z\otimes v_c+v_c\otimes z)\cdot z dx-\int_{\mathbb{R}^3}\text{div}(w\otimes w)\cdot z dx-\int_{\mathbb{R}^3} \nabla \pi_2 \cdot z dx.\label{4.7-0}
\end{eqnarray}
By similar argument as (\ref{aaaaa}), using integration by parts, we have
\begin{equation}
\begin{aligned}
-\int_{\mathbb{R}^3}\text{div}(z\otimes v_c+v_c\otimes z)\cdot  z dx\leq 2K_c\|\nabla z\|_{L^2}^2.
\end{aligned}
\end{equation}
By integration by parts, H\"{o}lder's inequality and Cauchy inequality, we have
\begin{equation}
-\int_{\mathbb{R}^3}\text{div}(w\otimes w)\cdot z dx=\int_{\mathbb{R}^3}(w\otimes w)\nabla z dx
\leq \|w\otimes w\|_{L^2}\|\nabla z \|_{L^2}
\leq C\|w\|_{L^4}^4+\frac {1}{10}\|\nabla z \|^2_{L^2}.
\end{equation}
For the pressure term, using integration by parts and div$z=0$, we get
\begin{equation}
-\int_{\mathbb{R}^3} \nabla \pi_2 \cdot z dx=0.\label{4.9-0}
 \end{equation}
 Then, from (\ref{4.7-0})-(\ref{4.9-0}), we have
\begin{equation} \label{4.15}
 \frac12\left\|z\right\|_{L_T^{\infty}L_x^2}^2+(\frac9{10}-2K_c)\left\|\nabla z\right\|_{L_T^2L_x^2}^2\leq C \|w\|_{L_T^{4}L_x^4}^4.
\end{equation}
Since the $L^p$ mild solution  $w\in L_T^{\infty}L_x^p\cap L_T^{\frac{4p}{3}}L_x^{2p}$,  $p\in[3,4]$, by interpolation theory, we have
$
w\in L_T^{\frac{8p}{3(4-p)}}L_x^4
$, and  $z\in  L^\infty([0,T];L_{x}^2)\cap L_T^2 (\dot{H}_{x}^1)$.
Combining with similar argument as (\ref{con_1}), we have $z\in  C([0,T];L_{x}^2)\cap L_T^2 (\dot{H}_{x}^1)$. Then, $w\in  C_w([0,T];L_{x}^2)\cap L_T^2 (\dot{H}_{x}^1)$,  one can easily prove that $w$ is a $L^2$-weak solution of system (\ref{PNS}) on $[0,T]$, and omit the details.

When $p>4$, the $L^p$ mild solution  $w\in L_T^{\infty}L_x^p\cap L_T^{\frac{4p}{3}}L_x^{2p}$, from Lemma \ref{p_lem3}, we could obtain that
    \begin{equation}\label{4.10}
\|z\|_{C_tL^\frac{p}{2}\cap L^{\frac{2p}{3}}_tL^{p}_x}\leq  C \|w\|_{L^\frac{2p}{p-3}_tL^{p}_x}\|w\|_{L^{\frac{2p}{3}}_tL^{p}_x},
\end{equation}
and $z\in  C([0,T];L_{x}^{\frac{p}{2}})\cap L^{\frac{2p}{3}}_tL^{p}_x$. Combing $a\in L_T^{\infty}L_x^p\cap L_T^{\frac{4p}{3}}L_x^{2p}\cap  C([0,T];L_{x}^2)\cap L_T^2 (\dot{H}_{x}^1)$, we have that $w\in  C([0,T];L_{x}^{\frac{p}{2}})\cap L^{\frac{2p}{3}}_tL^{p}_x$. By the induction, we can get
    $w\in  C([0,T];L_{x}^{\frac{p}{2^K}})\cap L^{\frac{p2^{2-K}}{3}}_tL^{ p2^{1-K}}_x$, for some $K\in \mathbb{Z}^+$ such that $2<\frac{p}{2^K}\leq 4$.
From the argument in (\ref{4.15}), we have that $z\in  C([0,T];L_{x}^2)\cap L_T^2 (\dot{H}_{x}^1)$.

Therefore, the $L^p$ mild solution $w\in  C_w([0,T];L_{x}^2)\cap L_T^2 (\dot{H}_{x}^1)$,  one can easily prove that $w$ is a $L^2$-weak solution of system (\ref{PNS}) on $[0,T]$, and omit the details.
  {\hfill
$\square$\medskip}

Based on the proof of Proposition \ref{prop-s-w-uniqueness}, we have the following results in global time. For simplicity, we omit the detailed proof.

\begin{corollary}\label{cor2}
For $p\geq 3,$ $T>0,$ let  ${c}_p$ and $\varepsilon_{0}$ be as in Theorem \ref{p>3 result}, $|c|>{c}_p$. For $w_0\in L_{\sigma}^p(\mathbb{R}^3)\cap L_{\sigma}^2(\mathbb{R}^3)$ and $\left\|{w}_{0}\right\|_{L^{3}\left(\mathbb{R}^{3}\right)}<\varepsilon_{0} ,$ let $w$ be a global  $L^p$ mild solution of system (\ref{PNS}). Then $w$ is a global $L^2$-weak solution of system (\ref{PNS}).
\end{corollary}

Combining with  Theorem \ref{p>3 result} and Corollary \ref{cor2}, we deduce Corollary \ref{corollary}.

\section{ Global $L^2+L^3$ weak solution}\label{Sec3}

 In this section, we will  illustrate  Theorem \ref{w-global-weak-existence}, i.e. we will give the global existence of $L^2+L^3$ weak solution to system (\ref{PNS}).

 Note that when we consider  the existence of weak solution to the Navier-Stokes system, there are essentially two methods: the energy method and the perturbation theory. The energy method gives the global existence for any initial data $v_0\in L_{\sigma}^2(\mathbb{R}^3)$. We cannot use this method since space $L^2$ doesn't contain space $L^3$ in the whole space $\mathbb{R}^3$. In the perturbation theory, by means of contraction mapping theorem, there exists a unique global weak solution to the Navier-Stokes system for small initial data $v_0\in L_{\sigma}^3(\mathbb{R}^3)$. Both methods cannot give direct results on the global existence for arbitrary $v_0\in L_{\sigma}^3(\mathbb{R}^3)$.

 Hence, many authors have developed various approaches to adapt the theory of the weak solutions so that it could allow $v_0\in L_{\sigma}^3(\mathbb{R}^3)$. Calder{\'o}n \cite{Cal} raised a method such that a $L_{\sigma}^3(\mathbb{R}^3)$ initial data $v_0$ can be decomposed as
\begin{equation}
v_0=v_0^1+v_0^2,
\end{equation}
where $v_0^1$ is small in $L_{\sigma}^3(\mathbb{R}^3)$ and $v_0^2$ belongs to $L_{\sigma}^2\cap L_{\sigma}^3(\mathbb{R}^3)$. Because of the smallness, initial data $v_0^1$ generates a global smooth solution $v_1$ by perturbation theory. Then the equation for $v_2=v-v_1$ can be solved by energy method. Seregin and $\check {\mathrm S}$ver$\acute {\mathrm a}$k \cite{Se} used another method to obtain global weak solution for $v_0\in L_{\sigma}^3(\mathbb{R}^3)$. The main idea of \cite{Se} is as follows. Let $v_1$ be  solution of the linear version of the Navier-Stokes system, seek  solution $v$ of the Navier-Stokes system as $v=v_1+v_2$, write down the equation that $v_2$ satisfied, then get the property of $v$ by investigating $v_2$. It's a general idea that the correction term $v_2$ might be easier to deal with than the full solution $v$. Related work can be referred in \cite{Lem,Ler,Se}.

 Inspired by above methods, we will decompose initial data $w_0=v_{10}+v_{20}$ and investigate the global existence of solutions $w=v_1+v_2$ to system (\ref{PNS}).
 For $w_0\in L_{\sigma}^3(\mathbb{R}^3)$, we have the following decomposition
\begin{equation}
w_0=v_{10}+v_{20},
\end{equation}
with $\|v_{10}\|_{L^3}<\varepsilon_0$ and $v_{20}\in L_{\sigma}^2\cap L_{\sigma}^3(\mathbb{R}^3).$ Since $\|v_{10}\|_{L^3}\ll1,$ there exists a unique global $L^3$ mild solution $v_1$ to system (\ref{v1}) according to Theorem \ref{L^3 well-posedness}. Crucial part is the global existence of $v_2$.
Since $v_{20}\in L_{\sigma}^2\cap L_{\sigma}^3$, this is the standard reasoning based on the Galerkin method (cf. \cite{Ka} Proof of Theorem 2.7). We claim there exist a
global weak solution $v_2\in C_{w}\left([0, T]; L_{\sigma}^{2}\left(\mathbb{R}^{3}\right)\right) \cap L^{2}\left([0, T]; \dot{H}_{\sigma}^{1}\left(\mathbb{R}^{3}\right)\right)$ for each $T>0.$ According to Definition \ref{global 3}, there exists
a global $L^2+L^3$ weak solution to system (\ref{PNS}). Detailed proof of the global existence of $v_2$ can be seen below.

First, we will construct weak solutions $v_2$ to the system (\ref{v2}). This is the standard reasoning based on the Galerkin method (cf. \cite{Ka} Proof of Theorem 2.1). Since $H_{\sigma}^{1}\left(\mathbb{R}^{3}\right)$ is separable, there exists a sequence $\left\{g_{m}\right\}_{m=1}^{\infty}$ which is free and total in $H_{\sigma}^{1}\left(\mathbb{R}^{3}\right)$. For each $m=1,2,\ldots$ Define an approximate solution $w_{m}=\sum_{i=1}^{m} d_{i m}(t) g_{i},$ which satisfies the following system of ordinary differential equations
\begin{equation}\label{l}
\begin{aligned}
&\left<w_{m}^{\prime}(t), g_{j}\right>+\left<\nabla w_{m}(t), \nabla g_{j}\right>+\left<\left(w_{m}(t) \cdot \nabla\right) w_{m}(t), g_{j}\right>\\
&+\left<\left(w_{m}(t) \cdot \nabla\right) (v_{c}+v_1), g_{j}\right> +\left<\left((v_{c}+v_1) \cdot \nabla\right) w_{m}(t), g_{j}\right>=0 \text { for } j=1, \ldots, m,
\end{aligned}
\end{equation}
where the term corresponding to the pressure in (\ref{v2}) vanishes in (\ref{l}) because of   $\text{div} g_{j}=0$.

We will prove terms $\left<\left(w_{m}(t) \cdot \nabla\right) (v_{c}+v_1), g_{j}\right> $ and $\left<\left((v_{c}+v_1) \cdot \nabla\right) w_{m}(t), g_{j}\right>$ in (\ref{l}) are convergent.
By H\"{o}lder and Sobolev inequalities in the Lorentz $L^{p, q} \text { -spaces (see \cite{Ka})}$, we have
\begin{eqnarray}\label{e}
\left|\int_{\mathbb{R}^{3}} g_j (w_m\cdot \nabla) (v_c+v_1) \mathrm{d} x\right|
&\leq& C\|(v_c+v_1)w_m\|_{L^2}\|\nabla g_j\|_{L^2} \nonumber\\
&\leq& C\|(v_c+v_1)w_m\|_{L^{2,2}}\|\nabla g_j\|_{L^2} \nonumber\\
&\leq& C\|v_c+v_1\|_{L^{3, \infty}}\|w_m\|_{L^{6,2}}\|\nabla g_j\|_{L^2} \nonumber\\
&\leq& C\|v_c+v_1\|_{L^{3, \infty}}\|\nabla w_m\|_{L^2}\|\nabla g_j\|_{L^2}.
\end{eqnarray}
Similar estimate holds for $\left<\left((v_{c}+v_1) \cdot \nabla\right) w_{m}(t), g_{j}\right>$.
\begin{equation}\label{e_2}
\begin{aligned}
\left|\int_{\mathbb{R}^{3}} ((v_c+v_1)\cdot \nabla) w_m g_j \mathrm{d} x\right| \leq C\|v_c+v_1\|_{L^{3, \infty}}\|\nabla w_m\|_{L^2}\|\nabla g_j\|_{L^2}.
\end{aligned}
\end{equation}

The system  (\ref{l}) has a unique local solution $\left\{d_{i m}(t)\right\}_{i=1}^{m} .$ By $a$-$priori$ estimates of the sequence $\left\{w_{m}\right\}_{m=1}^{\infty}$ obtained below in (\ref{s}),  solution $d_{i m}(t)$ is global.

Multiplying equation (\ref{l}) by $d_{j m}$ and sum up  equations for $j=1,2, \ldots, m ,$ we have
\begin{equation}
\frac{1}{2} \frac{\mathrm{d}}{\mathrm{d} t}\left\|w_{m}(t)\right\|_{2}^{2}+\left\|\nabla w_{m}(t)\right\|_{2}^{2}+\left<\left(w_{m}(t) \cdot \nabla\right)(v_c+v_1) , w_{m}(t)\right>=0.
\end{equation}
Using inequality (\ref{e}) and integrating it from $0$ to $t$, we obtain
\begin{equation}\label{s}
\left\|w_{m}(t)\right\|_{2}^{2}+2\left(1-K \sup _{t>0}\|v_c+v_1\|_{L^{3,\infty}}\right) \int_{0}^{t}\left\|\nabla w_{m}(\tau)\right\|_{2}^{2} \mathrm{d} \tau \leq \left\|w_{0}\right\|_{2}^{2}.
\end{equation}
Since $|c|$ big enough such that $K \sup _{t>0}\|v_c+v_1\|_{L^{3,\infty}}<1$. Thus we obtain a subsequence, also denoted by $\left\{w_{m}\right\}_{m=1}^{\infty}$, converging to $v_2\in C_{w}\left([0, T]; L_{\sigma}^{2}\left(\mathbb{R}^{3}\right)\right) \cap L^{2}\left([0, T]; \dot{H}_{\sigma}^{1}\left(\mathbb{R}^{3}\right)\right).$
Now, repeating the classical reasoning from \cite{Ka}, we obtain the existence of a weak solution in the energy space $C_{w}\left([0, T]; L_{\sigma}^{2}\left(\mathbb{R}^{3}\right)\right) \cap L^{2}\left([0, T]; \dot{H}_{\sigma}^{1}\left(\mathbb{R}^{3}\right)\right)$ for all $T>0$ which satisfies strong energy inequality (\ref{s}).

Hence we get a global  $L^2+L^3$ weak solution $w$ of the form $w=v_1+v_2$.
Moreover, we have the asymptotic behavior of $v_2$ and omit the proof which can be referred in \cite{Ka}.
  {\hfill
$\square$\medskip}

\section{Proof of Theorem \ref{p>3 result}}\label{Sec4}
In this section, we will give the proof of Theorem \ref{p>3 result}. Our method is based on contraction mapping in Lemma \ref{lem1} and the following $a$-$prioi$  estimates in Lemmas \ref{p_alemm}, \ref{p_lem2} and \ref{p_lem3}.

\begin{lemma} \label{p_alemm}
 Let $p\in (1,\infty)$.  For  every $c$ satisfies (\ref{condition c_p1}) and (\ref{condition c_p2}), there exists a unique global-in-time  solution $a(x,t)\in C_tL_x^p\cap L^{\frac{4p}{3}}_tL^{2p}_x$ to system (\ref{a}) with initial data $w_0\in L_{\sigma}^p({\mathbb{R}^{3}}).$  Moreover,
 \begin{equation}\label{6.1}
\|a(\cdot,t)\|_{L^p}\leq \|a(\cdot,s)\|_{L^p},
 \end{equation}
 for any $0\leq t\leq s<\infty,$ and for $2\leq p<\infty$,
\begin{equation}\label{p_a-estimate}
\|a\|_{C_tL_x^p\cap L^{\frac{4p}{3}}_tL^{2p}_x}+\left\|\nabla\left(|a|^{\frac{p}{2}}\right)\right\|^{\frac2p}_{L_t^2L_x^2}\leq C\|w_0\|_{L^p},
\end{equation}
for a universal constant $C$.
\end{lemma}
\noindent\textbf{Proof of Lemma \ref{p_alemm}.}
 Approach is similar as the proof of Lemma \ref{alemm}. By classical approximation method, it is easy to get the global existence of solutions $a$. For simplicity, we omit the detailed proof and give $a$-$prioi$ estimate for $a.$ Suppose $a$ is sufficiently smooth, we multiply the equation ($\ref{a})_1$ by $|a|^{p-2}a$ and integrate it on $\mathbb{R}^3$, we have
\begin{equation}
\int_{\mathbb{R}^3}\partial_t a \cdot (|a|^{p-2}a)dx=\frac{1}{p}\frac{\mathrm{d}}{\mathrm{d} t}\|a(t)\|^p_{L^p},
\end{equation}
and
\begin{eqnarray}\label{000}
 \int_{\mathbb{R}^3}-\Delta a \cdot (|a|^{p-2}a)dx
&=&(p-2)\int_{\mathbb{R}^3}|a|^{p-4} \sum_i[(\partial_i a_l)a_l]^2 + \int_{\mathbb{R}^3} |\nabla a|^2|a|^{p-2}\nonumber\\
&=&\frac{4(p-2)}{p^2}\|\nabla(|a|^{\frac{p}{2}})\|^2_{L^2}+\||\nabla a||a|^{\frac{p-2}{2}}\|^2_{L^2}.
\end{eqnarray}
For $p\geq 2,$ we have
\begin{eqnarray}\label{11}
&&\frac{1}{p}\frac{\mathrm{d}}{\mathrm{d} t}\|a(t)\|^p_{L^p}+\frac{4(p-2)}{p^2}\|\nabla(|a|^{\frac{p}{2}})\|^2_{L^2}+\||\nabla a||a|^{\frac{p-2}{2}}\|^2_{L^2}\nonumber\\
&=&-\int_{\mathbb{R}^3}\text{div}(a\otimes v_c+v_c\otimes a)\cdot (|a|^{p-2}a) dx-\int_{\mathbb{R}^3} \nabla \pi\cdot (|a|^{p-2}a) dx.
\end{eqnarray}
By using  integration by parts, H\"{o}lder's inequality, Lemma \ref{|x|v_c} and the classical Hardy inequality in Lemma \ref{hardy}, we have
\begin{eqnarray}\label{22}
 -\int_{\mathbb{R}^3}\text{div}(a\otimes v_c+v_c\otimes a)\cdot (|a|^{p-2}a) dx
&=&\int_{\mathbb{R}^3}(a\otimes v_c+v_c\otimes a) \cdot \nabla(|a|^{p-2}a )dx\nonumber\\
&=& \int_{\mathbb{R}^3} a_i(v_c)_j\partial_i(|a|^{p-2}a_j)   dx\nonumber\\
&\leq& C\int_{\mathbb{R}^3}|\nabla(|a|^{\frac{p}{2}} )||a|^{\frac{p}{2}} |v_c|dx\nonumber\\
&\leq& C\left\||x| v_c\right\|_{L^{\infty}}\left\|\nabla(|a|^{\frac{p}{2}} )\right\|_{L^2}\left\|\frac{|a|^{\frac{p}{2}}}{|x|}\right\|_{L^2}\nonumber\\
&\leq& CK_c\left\|\nabla(|a|^{\frac{p}{2}} )\right\|_{L^2}^2.
\end{eqnarray}
 Similar as proof of (\ref{prop- rh2}), by  H\"{o}lder's inequality, Hardy inequality, Sobolev embedding and boundedness of the Riesz transforms on weighted $L^p$ spaces (Theorem 9.4.6 in \cite{Gra}), we have
\begin{eqnarray}\label{33}
 \int_{\mathbb{R}^3} \nabla \pi_1 \cdot (|a|^{p-2}a) dx
&\leq& C\||x|^{\frac{p-2}{p}}\left(a\otimes v_c+v_c\otimes a\right)\|_{L^{p}} \|\nabla (|a|^{\frac p2})\|_{L^{2}} \left\|\frac{|a|^{\frac p2-1}}{|x|^{\frac{p-2}{p}}}\right\|_{L^{\frac{2p}{p-2}}}\nonumber\\
&\leq& C\||x|v_c\|_{L^{\infty}} \left\|\frac{a}{|x|^{\frac2p}}\right\|_{L^p}\|\nabla (|a|^{\frac p2})\|_{L^{2}} \left\|\nabla(|a|^{\frac p2})\right\|_{L^{2}}^{\frac {p-2} {p} }\nonumber\\
&\leq& C K_c\left\|\frac{a}{|x|^{\frac2p}}\right\|_{L^p}\|\nabla (|a|^{\frac{p}{2}})\|_{L^2}^{\frac{2p-2}{p}}\nonumber\\
&\leq& C K_c \|\nabla (|a|^{\frac{p}{2}})\|^2_{L^2}.
\end{eqnarray}
Combining (\ref{11})-(\ref{33}), we deduce
\begin{equation}\label{a_p}
\frac 1p\frac{\mathrm{d}}{\mathrm{d} t}\|a(t)\|^p_{L^p}+\frac{4(p-2)}{p^2}\|\nabla(|a|^{\frac{p}{2}})\|^2_{L^2}\leq C K_c \|\nabla (|a|^{\frac{p}{2}})\|^2_{L^2}.
\end{equation}
By assumption, we can guarantee
\begin{equation}\label{condition c_p1}
\frac{4(p-2)}{p^2}-CK_c>0,
\end{equation}
for a constant $C.$
 Hence
\begin{equation}
\frac{\mathrm{d}}{\mathrm{d} t}\|a(t)\|^p_{L^p}+C\|\nabla(|a|^{\frac{p}{2}})\|^2_{L^2}\leq 0,
\end{equation}
for a positive constant $C.$
Therefore,  (\ref{6.1}) holds and we have
\begin{equation}
\sup_t\|a(t)\|^p_{L^p}+C\|\nabla(|a|^{\frac{p}{2}})\|^2_{L_t^2L^2_x}\leq \|w_0\|_{L^p}.
\end{equation}
By interpolation theory, we deduce (\ref{p_a-estimate}).

By Lemma \ref{sec3-lem1}, we obtain
\begin{equation}
\sum_i[(\partial_i a_l)a_l]^2 \leq |a|^2|\nabla a|^2.
\end{equation}
For $1< p<2,$  we have
\begin{eqnarray*}
&&\frac{1}{p}\frac{\mathrm{d}}{\mathrm{d} t}\|a(t)\|^p_{L^p}+(p-2)\int_{\mathbb{R}^3}|a|^{p-4} \sum_i[(\partial_i a_l)a_l]^2 + \int_{\mathbb{R}^3} |\nabla a|^2|a|^{p-2}\\
&=&-\int_{\mathbb{R}^3}\text{div}(a\otimes v_c+v_c\otimes a)\cdot (|a|^{p-2}a) dx-\int_{\mathbb{R}^3} \nabla \pi\cdot (|a|^{p-2}a) dx.
\end{eqnarray*}
Thanks to Lemma \ref{sec3-lem1}, there holds
\begin{eqnarray*}
&&\frac{1}{p}\frac{\mathrm{d}}{\mathrm{d} t}\|a(t)\|^p_{L^p}+(p-1)\int_{\mathbb{R}^3} |\nabla a|^2|a|^{p-2}\\
&\leq&\Big|\int_{\mathbb{R}^3}\text{div}(a\otimes v_c+v_c\otimes a)\cdot (|a|^{p-2}a) dx\Big|+\Big|\int_{\mathbb{R}^3} \nabla \pi\cdot (|a|^{p-2}a) dx\Big|.
\end{eqnarray*}
Thanks to (\ref{22}) and (\ref{33}), there holds
\begin{equation}
\begin{aligned}
&\frac{1}{p}\frac{\mathrm{d}}{\mathrm{d} t}\|a(t)\|^p_{L^p}+(p-1)\int_{\mathbb{R}^3} |\nabla a|^2|a|^{p-2}dx\leq C K_c \|\nabla (|a|^{\frac{p}{2}})\|^2_{L^2}.
\end{aligned}
\end{equation}
Moreover,
\begin{equation}
\begin{aligned}
C K_c \|\nabla (|a|^{\frac{p}{2}})\|^2_{L^2}&= C K_c \frac{p^2}{4}\int_{\mathbb{R}^3}|a|^{p-4} \sum_i[(\partial_i a_l)a_l]^2 dx\\
&\leq C(p)K_c \int_{\mathbb{R}^3}|a|^{p-2} |\nabla a|^2dx,
\end{aligned}
\end{equation}
where $C(p)$ is a constant depending on $p.$
By assumption, we can guarantee
\begin{equation}\label{condition c_p2}
p-1-C(p)K_c>0.
\end{equation}
Therefore
\begin{equation}
\begin{aligned}
&\frac{\mathrm{d}}{\mathrm{d} t}\|a(t)\|^p_{L^p}+C \int_{\mathbb{R}^3}|a|^{p-2} |\nabla a|^2dx\leq 0,
\end{aligned}
\end{equation}
for a positive constant $C.$
Hence we deduce (\ref{6.1}).
 {\hfill
$\square$\medskip}

To get $a$-$prioi$  estimate of $z$,  in which the crucial estimate is as follows:
\begin{eqnarray*}
 -\int_{\mathbb{R}^3}\text{div}(w_1\otimes w_2) \cdot (|z|^{p-2}z) dx
&=&\int_{\mathbb{R}^3}(w_1\otimes w_2)\cdot  \nabla(|z|^{p-2}z )dx\\
&\leq& C\int_{\mathbb{R}^3}(w_1\otimes w_2)\cdot
\nabla\left(|z|^{\frac p2} \right)|z|^{\frac p2-1} dx\\
&\leq& C\left\|\nabla\left(|z|^{\frac{p}{2}} \right)\right\|_{L^2}\left\||z|^{\frac{p}{2}-1}\right\|_{L^{\frac{2p}{p-2}}}\left\|w_1\otimes w_2\right\|_{L^p}\\
&\leq& C\left\|\nabla\left(|z|^{\frac{p}{2}} \right)\right\|_{L^2}\left\|z\right\|_{L^p}^{\frac{p}{2}-1}\left\|w_1\otimes w_2\right\|_{L^p}\\
&\leq& \varepsilon \left\|\nabla\left(|z|^{\frac{p}{2}} \right)\right\|^2_{L^2}+C(\varepsilon)\left\|z\right\|_{L^p}^{p-2}\left\|w_1\otimes w_2\right\|^2_{L^p}.\\
\end{eqnarray*}
We have
\begin{eqnarray*}
\sup_t \|z(t)\|_{L_x^p}+\left\|\nabla(|z|^{\frac{p}{2}})\right\|^{\frac{2}{p}}_{L_t^2L_x^2} \leq  C \|w_1\|_{L_t^{4}L_x^{2p}}\|w_2\|_{L_t^{4}L_x^{2p}}
&\leq&C T^{\frac{p-3}{2p}} \|w_1\|_{ L^{\frac{4p}{3}}_t L^{2p}_x}\|w_2\|_{L^{\frac{4p}{3}}_t L^{2p}_x}.
\end{eqnarray*}
Hence, we have the following $a$-$prioi$  estimate and more detailed proof can be referred in the proof of  Lemma \ref{lem2}.

\begin{lemma}\label{p_lem2}
 Let $p\in (3,\infty)$. Assume that ${c}_p$ is as in Theorem \ref{p>3 result}. For every $|c|>{c}_p$, there exists a  $L^p$ mild solution $z(x,t)$ on $[0,T]$ to system (\ref{N}) with $w_1, w_2\in L_t^{\frac{4p}{3}}([0,T];L^{2p}(\mathbb{R}^3))$, satisfying
\begin{eqnarray}\label{p_z-estimate}
\|z\|_{C_TL^p\cap L^{\frac{4p}{3}}_TL^{2p}_x}+\|\nabla(|z|^{\frac p2})\|_{L_t^2L_x^2}^{\frac2p}\leq  C T^{\frac{p-3}{2p}} \|w_1\|_{L^{\frac{4p}{3}}_TL^{2p}_x}\|w_2\|_{L^{\frac{4p}{3}}_TL^{2p}_x},
\end{eqnarray}
for a constant $C$.
\end{lemma}
When initial data $w_0\in L_{\sigma}^p\cap L_{\sigma}^3$ and $\|w_0\|_{L^3}<\varepsilon_0$, we have $w\in C_tL^3\cap L^{4}_tL^6_x$, $\nabla\left(|w|^{\frac{p}{2}}\right)\in L_t^2L_x^2$ according Theorem \ref{L^3 well-posedness}. By interpolation, we have $w\in L^{\frac{4p}{2p-3}}_tL^{2p}_x.$ Proof is very similar as the proof of Lemma \ref{lem2}, in which the crucial estimate is as follows:
\begin{equation}
\|w_1\otimes w_2\|_{L_t^{2}L_x^p}\leq  \|w_1\|_{L^{\frac{4p}{2p-3}}_tL^{2p}_x}\|w_2\|_{L^{\frac{4p}{3}}_tL^{2p}_x}.
\end{equation}
Hence we have the following $a$-$prioi$  estimate and more detailed proof can be referred in the proof of  Lemma \ref{lem2}.
\begin{lemma}\label{p_lem3}
 Let $p\in (3,\infty)$. Assume that ${c}_p$ is as in Theorem \ref{p>3 result}. For
 every $|c|>{c}_p$, there exists a global-in-time $L^p$ mild solution $z(x,t)$ to
 system (\ref{N}) with $w_1\in L_t^{\frac{4p}{2p-3}}([0,\infty);L_x^{2p}(\mathbb{R}^3))$ and $w_2\in L_t^{\frac{4p}{3}}([0,\infty);L_x^{2p}(\mathbb{R}^3))$, satisfying
\begin{equation}\label{p_z-estimate2}
\|z\|_{C_tL^p\cap L^{\frac{4p}{3}}_tL^{2p}_x}+ \|\nabla(|z|^{\frac p2})\|_{L_t^2L_x^2}^{\frac2p}\leq  C \|w_1\|_{L^\frac{4p}{2p-3}_tL^{2p}_x}\|w_2\|_{L^{\frac{4p}{3}}_tL^{2p}_x},
\end{equation}
for a constant $C$.
\end{lemma}
\begin{remark}  \label{p_existence}
By classical method, it is easy to get the existence of solutions $a\in C([0,\infty);L^p(\mathbb{R}^3))$, $\nabla (|a|^{\frac p2})\in L^2([0,\infty);L^2(\mathbb{R}^3))$ to system (\ref{a}) satisfying (\ref{p_a-estimate}) and solutions $z\in C([0,\infty);L^p(\mathbb{R}^3))$, $\nabla (|z|^{\frac p2})\in L^2([0,\infty);L^2(\mathbb{R}^3))$ to system (\ref{N}). For simplicity, we omit the detailed proof.
\end{remark}

\noindent\textbf{Proof of Theorem \ref{p>3 result}.}
 For constant $|c|>{c}_p$ where  ${c}_p$  depends only on $p$,
according to Lemma  \ref{p_alemm}, we have
\begin{equation}
\|a(t)\|_{ L^{\frac{4p}{3}}_tL^{2p}_x}\leq C \|w_0\|_{L^p}.
\end{equation}
Applying Lemma  \ref{p_lem2} with $(w_1,w_2)=(w,w)$, we have
\begin{equation}
\|N\|\leq CT^{\frac{p-3}{2p}}.
\end{equation}
Using Lemma \ref{lem1} with $E= L^{\frac{4p}{3}}_TL^{2p}_x$, there exists $T>0$ and a unique solution $w\in L^{\frac{4p}{3}}_TL^{2p}_x$ on $[0,T].$

Then we will prove the global existence of $w$ with initial data $w_0\in L^p_{\sigma}\cap L_{\sigma}^3$ and $\|w_0\|_{L^3}<\varepsilon_0$. Since $\|w_0\|_{L^3}<\varepsilon_0$, according to Theorem \ref{L^3 well-posedness}, there exists a global unique solution $w\in C_tL^{3}_x\cap L^{4}_tL^6_x,$ $\nabla (|w|^{\frac 32})\in L_t^2 L_x^2,$ and $\|w\|_{C_tL^{3}_x\cap L^{4}_tL^6_x}+\|\nabla |w|^{\frac 23}\|^{\frac 23}_{L_t^2L_x^2}\leq C\|w_0\|_{L^3}.$  By interpolation, $w\in C_t{L_x^3}$ and  $\nabla (|w|^{\frac 32})\in L_t^2 L_x^2$  deduce $w\in L^{\frac{4p}{2p-3}}_tL^{2p}_x$.
Hence,
\begin{equation}\label{3-p}
\|w\|_{L^{\frac{4p}{2p-3}}_tL^{2p}_x}\leq C\|w_0\|_{L^3}< C\varepsilon_0 .
\end{equation}
Thanks to (\ref{p_a-estimate}) and (\ref{p_z-estimate2}), we have
\begin{equation}\label{w_p}
\|w\|_{C_tL^{p}_x\cap L^{\frac{4p}{3}}_tL^{2p}_x}\leq C\|w_0\|_{L^p}+C \|w\|_{L^{\frac{4p}{2p-3}}_tL^{2p}_x}\|w\|_{C_tL^{p}_x\cap L^{\frac{4p}{3}}_tL^{2p}_x}.
\end{equation}
Combining with (\ref{3-p}) and interpolation theory, we deduce (\ref{p>3-1}).
 {\hfill
$\square$\medskip}

\section{Proof of Theorems \ref{continuty} }\label{proof of {continuty} }

In this section, we give the detailed proof of Theorems \ref{continuty}.

\noindent\textbf{Proof of Theorem \ref{continuty}.}\\
Setting $Z=u-v$, we have
\begin{equation}\label{u-v}
\begin{cases}
 Z_{t}-\Delta Z + \text{div}(-Z \otimes Z+Z \otimes u+u\otimes Z )+ (Z  \cdot \nabla) v_{c}+\left(v_{c} \cdot \nabla\right) z+\nabla \pi_{z}=0,\\
 \nabla\cdot Z=0,\\
 Z(x, 0) =Z_0.
\end{cases}
\end{equation}
By the Duhamel principle, we can write solution $z$ into an integral formulation
\begin{equation}
Z(x,t)=e^{-t\mathcal{L}}u_0-\int_0^t e^{-(t-s)\mathcal{L}}\mathbb{P}\text{div}(-Z \otimes Z+Z \otimes u+u\otimes Z )ds.
\end{equation}
By contraction mapping theorem, it's easy to give the existence of solution. Next, we only give  $a$-$prioi$ estimate.

When  $p=3,$
by Lemma \ref{alemm} and method in Lemma \ref{lem2}, we have
\begin{equation}\label{eq1}
 \|Z\|_{C_TL_x^3\cap L^4_TL_x^6 }+\left\|\nabla\left(|Z|^{\frac{3}{2}}\right)\right\|^{\frac23}_{L_T^2L_x^2}
\leq   C_1\|Z_0\|_{L^3}+C_2\|Z\|^2_{L^4_TL_x^6}+C_2\left(\int_0^T (\|Z\|_{L^6}\|u\|_{L^6})^2 dt\right)^{\frac12}.
\end{equation}
By interpolation inequality, H\"{o}lder's inequality and Young's inequality, we have
\begin{eqnarray*}
\left(\int_0^T (\|Z\|_{L^6}\|u\|_{L^6})^2 dt\right)^{\frac12}
&\leq& \left(\int_0^T \left(\|Z\|^{\frac14}_{L^3}\|Z\|^{\frac34}_{L^9}\|u\|_{L^6}\right)^2 dt\right)^{\frac12}\\
&\leq& \left(\int_0^T \|Z\|^{\frac12}_{L^3}\|Z\|^{\frac 32}_{L^9}\|u\|_{L^6}^2 dt\right)^{\frac12}\\
&\leq& \|Z\|^{\frac 34}_{L_T^3L_x^9} \left(\int_0^T \|Z\|_{L^3}\|u\|_{L^6}^4 dt\right)^{\frac14}\\
&\leq& \varepsilon   \|Z\|_{L_T^3L_x^9} +\frac{27}{256 \varepsilon^3} \int_0^T \|Z\|_{L^3}\|u\|_{L^6}^4 dt.
\end{eqnarray*}
Combining with (\ref{eq1}), we obtain
\begin{eqnarray*}
 \|Z\|_{C_TL_x^3\cap L^4_TL_x^6 }+\left\|\nabla\left(|Z|^{\frac{3}{2}}\right)\right\|^{\frac23}_{L_T^2L_x^2}
&\leq&  C\|Z_0\|_{L^3}+C\|Z\|^2_{L^4_TL_x^6}+C \varepsilon \|Z\|_{L_T^3L_x^9} +\frac{C}{ \varepsilon^3} \int_0^T \|Z\|_{L^3}\|u\|_{L^6}^4 dt.
\end{eqnarray*}
 By Sobolev embedding $\dot{H}^{1}(\mathbb{R}^{3})\hookrightarrow L^6(\mathbb{R}^{3})$, we have
\begin{eqnarray*}
&&\|Z\|_{C_TL_x^3\cap L^4_TL_x^6 }+\left\|\nabla\left(|Z|^{\frac{3}{2}}\right)\right\|^{\frac23}_{L_T^2L_x^2}\\
&\leq&  C\|Z_0\|_{L^3}+C\|Z\|^2_{L^4_TL_x^6}+C \varepsilon\left\|\nabla\left(|Z|^{\frac{3}{2}}\right)\right\|^{\frac23}_{L_T^2L_x^2} +\frac{C}{ \varepsilon^3} \int_0^T \|Z\|_{L^3}\|u\|_{L^6}^4 dt.
\end{eqnarray*}
 Taking $C \varepsilon=\frac12$, there holds
\begin{eqnarray}\label{*}
 \|Z\|_{C_TL_x^3\cap L^4_TL_x^6 }+\left\|\nabla\left(|Z|^{\frac{3}{2}}\right)\right\|^{\frac23}_{L_T^2L_x^2}
&\leq&  C\|Z_0\|_{L^3}+C\|Z\|^2_{L^4_TL_x^6} +C \int_0^T \|Z\|_{L^3}\|u\|_{L^6}^4 dt,
\end{eqnarray}
for a positive constant $C.$
According to Gronwall's inequality, there holds
\begin{eqnarray}
 \|Z\|_{C_TL_x^3\cap L^4_TL_x^6 }+\left\|\nabla\left(|Z|^{\frac{3}{2}}\right)\right\|^{\frac23}_{L_T^2L_x^2}
&\leq&  C(\|Z_0\|_{L^3}+\|Z\|^2_{L^4_TL_x^6})e^{C \int_0^T \|u\|_{L^6}^4 dt}.
\end{eqnarray}
When $\|Z_0\|_{L^3}\leq (4C^2e^{2C \int_0^T \|u\|_{L^6}^4 dt})^{-1},$ by continuity method, we have
\begin{eqnarray}\label{Z-crucial}
\|Z\|_{C_{T}L_x^3\cap L^4_{T}L_x^6 }+\left\|\nabla\left(|Z|^{\frac{3}{2}}\right)\right\|^{\frac23}_{L_{T}^2L_x^2}\leq  2C\|Z_0\|_{L^3}e^{C\int_0^{T} \|u\|_{L^6}^4 dt}.
\end{eqnarray}
Therefore,  (\ref{sol}) holds with $p=3$.

When  $p>3,$
by Lemma \ref{p_alemm} and method in Lemma \ref{p_lem2}, we have
\begin{eqnarray}\label{eq2}
&&\|Z\|_{C_TL_x^p\cap L^{\frac{4p}{3}}_TL_x^{2p} }+\left\|\nabla\left(|Z|^{\frac{p}{2}}\right)\right\|^{\frac2p}_{L_T^2L_x^2}\\
&\leq&   C\|Z_0\|_{L^p}+C \left(\int_0^T (\|Z\|_{L_x^{2p}}\|Z\|_{L_x^{2p}})^2 dt\right)^{\frac12}+C \left(\int_0^T (\|Z\|_{L_x^{2p}}\|u\|_{L_x^{2p}})^2 dt\right)^{\frac12}.\nonumber
\end{eqnarray}
By interpolation inequality, H\"{o}lder's inequality and Young's inequality, we have
\begin{eqnarray*}
\left(\int_0^T \left(\|Z\|_{L_x^{2p}}\|u\|_{L_x^{2p}}\right)^2 dt\right)^{\frac12}
&\leq&  \left(\int_0^T \left(\|Z\|^{\frac 34}_{L_x^{3p}}\|Z\|^{\frac 14}_{L_x^p}\|u\|_{L_x^{2p}}\right)^2 dt\right)^{\frac{1}{2}}\\
&\leq&  \left\|\|Z\|^{\frac 34}_{L_x^{3p}}\right\|_{L_T^{\frac{4p}{3}}} \left\|\|Z\|^{\frac 14}_{L_x^{p}}\|u\|_{L_x^{2p}}\right\|_{L_T^{\frac{4p}{2p-3}}}\\
&\leq&  \|Z\|^{\frac 34}_{L_T^{p}L_x^{3p}} \left\|\|Z\|^{\frac 14}_{L_x^{p}}\|u\|_{L_x^{2p}}\right\|_{L_T^{\frac{4p}{2p-3}}}\\
&\leq& \varepsilon\|Z\|_{L_T^pL_x^{3p}} +\frac{27}{256 \varepsilon^3}\left\|\|Z\|^{\frac 14}_{L_x^{p}}\|u\|_{L_x^{2p}}\right\|^4_{L_T^{\frac{4p}{2p-3}}}.
\end{eqnarray*}
Combining with (\ref{eq2}), we obtain
\begin{eqnarray*}
&&\|Z\|_{C_TL_x^p\cap L^{\frac{4p}{3}}_TL_x^{2p} }+\left\|\nabla\left(|Z|^{\frac{p}{2}}\right)\right\|^{\frac2p}_{L_T^2L_x^2}\\
&\leq&   C\|Z_0\|_{L^p}+C \left(\int_0^T (\|Z\|_{L_x^{2p}}\|Z\|_{L_x^{2p}})^2 dt\right)^{\frac12}
 + C\varepsilon\|Z\|_{L_T^pL_x^{3p}} +\frac{C}{ \varepsilon^3}\left\|\|Z\|^{\frac 14}_{L_x^{p}}\|u\|_{L_x^{2p}}\right\|^4_{L_T^{\frac{4p}{2p-3}}}.
\end{eqnarray*}
 By Sobolev embedding $\dot{H}^{1}(\mathbb{R}^{3})\hookrightarrow L^6(\mathbb{R}^{3})$, we have
\begin{eqnarray*}
&&\|Z\|_{C_TL_x^p\cap L^{\frac{4p}{3}}_TL_x^{2p} }+\left\|\nabla\left(|Z|^{\frac{p}{2}}\right)\right\|^{\frac2p}_{L_T^2L_x^2}\\
&\leq&   C\|Z_0\|_{L^p}+C \left(\int_0^T (\|Z\|_{L_x^{2p}}\|Z\|_{L_x^{2p}})^2 dt\right)^{\frac12}
 + C\varepsilon\left\|\nabla\left(|Z|^{\frac{p}{2}}\right)\right\|^{\frac2p}_{L_T^2L_x^2} +\frac{C}{ \varepsilon^3}\left\|\|Z\|^{\frac 14}_{L_x^{p}}\|u\|_{L_x^{2p}}\right\|^4_{L_T^{\frac{4p}{2p-3}}}.
\end{eqnarray*}
Taking $C\varepsilon=\frac 12$, there holds
\begin{eqnarray}\label{**}
&&\|Z\|_{C_TL_x^p\cap L^{\frac{4p}{3}}_TL_x^{2p} }+\left\|\nabla\left(|Z|^{\frac{p}{2}}\right)\right\|^{\frac2p}_{L_T^2L_x^2}\nonumber\\
&\leq&   C\|Z_0\|_{L^p}+C \left(\int_0^T (\|Z\|_{L_x^{2p}}\|Z\|_{L_x^{2p}})^2 dt\right)^{\frac12}+C\left\|\|Z\|^{\frac 14}_{L_x^{p}}\|u\|_{L_x^{2p}}\right\|^4_{L_T^{\frac{4p}{2p-3}}}\nonumber\\
&\leq&   C\|Z_0\|_{L^p}+CT^{\frac{p-3}{2p}}\|Z\|^2_{L^{\frac{4p}{3}}_TL_x^{2p} }+C\|Z\|_{L_T^{\infty}L_x^{p}}\|u\|^4_{L_T^{\frac{4p}{2p-3}}L_x^{2p}},
\end{eqnarray}
for a positive constant $C.$
Thanks to Gronwall's inequality, there holds
$$
 \|Z\|_{C_TL_x^p\cap L^{\frac{4p}{3}}_TL_x^{2p} }+\left\|\nabla\left(|Z|^{\frac{p}{2}}\right)\right\|^{\frac2p}_{L_T^2L_x^2}
 \leq   C\left(\|Z_0\|_{L^p}+T^{\frac{p-3}{2p}}\|Z\|^2_{L^{\frac{4p}{3}}_TL_x^{2p} }\right)\exp\left\{{C\|u\|^4_{L_T^{\frac{4p}{2p-3}}L_x^{2p}}}\right\}.
$$
By the continuity method, when $4C^2T^{\frac{p-3}{2p}}\|Z_0\|_{L^p}\exp \{{2C\|u\|^4_{L_T^{\frac{4p}{2p-3}}L_x^{2p}}} \}<1$, we deduce
\begin{equation}\label{qq}
\begin{aligned}
\|Z\|_{C_{T}L_x^p\cap L^{\frac{4p}{3}}_{T}L_x^{2p} }+\left\|\nabla\left(|Z|^{\frac{p}{2}}\right)\right\|^{\frac2p}_{L_{T}^2L_x^2}\leq  2C\|Z_0\|_{L^p}\exp\left\{{C\|u\|^4_{L_T^{\frac{4p}{2p-3}}L_x^{2p}}}\right\}.
\end{aligned}
\end{equation}
When $\|Z_0\|_{L^p}\rightarrow0$, (\ref{qq}) implies  (\ref{sol})  with $p>3$.
 {\hfill
$\square$\medskip}

\section*{Acknowledgment}
The work of the first named author is partially supported by NSF  grant  DMS-1501004 and DMS-2000261. The work of the third named author is partially supported by the NSFC grant 11771389,  11931010 and 11621101.


%
%


\begin{thebibliography}{}
%
%

\bibitem{Ba} 
  \newblock H. Bahouri, J.-Y. Chemin, R. Danchin,
     \newblock \emph{Fourier Analysis and Nonlinear Partial Differential Equations},
     \newblock   Grundlehren Math. Wiss., vol. 343, Springer-Verlag, Berlin, Heidelberg, 2011.

  \bibitem{Bas}
\newblock A. Basson,
\newblock \emph{Solutions spatialement homog$\grave{e}$nes adapt$\acute{e}$es des
 $\acute{e}$quations de {N}avier-{S}tokes},
\newblock Thesis. University of Evry., 2006.



\bibitem{Bor1} 
\newblock W. Borchers and T. Miyakawa,
\newblock $L^2$ decay for {N}avier-{S}tokes flows in unbounded domains, with application to exterior stationary flows,
\newblock  \emph{Arch. Ration. Mech. Anal.} \textbf{118} (1992), 273-295.

\bibitem{Bo} 
\newblock W. Borchers, T. Miyakawa,
\newblock On stability of exterior stationary  {N}avier-{S}tokes flows,
\newblock  \emph{Acta Math.} \textbf{174} (1995), 311-382.




\bibitem{Caf} 
\newblock L. Caffarelli, R. Kohn, L. Nirenberg,
\newblock Partial regularity of suitable weak solutions of the Navier-Stokes equations,
\newblock  \emph{Comm. Pure Appl. Math.} \textbf{35} (6) (1982), 771-831.



\bibitem{Cal} 
\newblock C.P. Calder{\'o}n,
\newblock Existence of weak solutions for the Navier-Stokes equations with initial data in $L^p$,
\newblock  \emph{Trans. Am. Math. Soc.} \textbf{318} (1) (1990), 179-200.




\bibitem{Can0} 
  \newblock M. Cannone,
     \newblock \emph{Ondelettes, paraproduits et {N}avier-{S}tokes,}
     \newblock Diderot Editeur, Paris, 1995.






\bibitem{Can1} 
  \newblock M. Cannone,
     \newblock \emph{Harmonic analysis tools for solving the incompressible {N}avier-{S}tokes equations},
     \newblock in: Handbook of Mathematical Fluid Dynamics, vol. {III}, North-Holland, Amsterdam, 2004, pp. 161-244.




\bibitem{Can}
\newblock M. Cannone, G. Karch,
\newblock Smooth or singular solutions to the {N}avier-{S}tokes system?
\newblock  \emph{J. Differential Equations} \textbf{197} (2004), 247-274.






\bibitem{Ca}
\newblock E.A. Carlen, M. Loss,
\newblock Optimal smoothing and decay estimates for viscously damped conservation laws, with applications to the 2-D Navier-Stokes equation,
\newblock  \emph{Duke Math. J.} \textbf{81} (1995), 135-157.


\bibitem{Gra} 
  \newblock L. Grafakos,
     \newblock \emph{Modern Fourier analysis},
     \newblock  Graduate Texts in Mathematics, vol. 250, Springer, New York, 2009.






 \bibitem{Har1}
\newblock G.H. Hardy,
\newblock Note on a theorem of {H}ilbert,
 \newblock  \emph{Math. Zeit.} \textbf{6} (1920), 314-317.

 \bibitem{Har2}
\newblock G.H. Hardy,
\newblock An inequality between integrals,
 \newblock \emph{Messenger of Mathematics,} \textbf{54} (1925), 150-156.
 \bibitem{Iwa}
\newblock T. Iwaniec, G. Martin,
\newblock \emph{Riesz transforms and related singular integrals},
 \newblock  \emph{J. Reine Angew. Math.} \textbf{473} (1996), 25-57.

 \bibitem{Jia}
\newblock H. Jia, V. \v{S}ver\'{a}k,
\newblock Are the incompressible 3D {N}avier-{S}tokes equations locally
              ill-posed in the natural energy space?,
 \newblock  \emph{J. Funct. Anal.} \textbf{268} (12) (2015), 3734-3766.


\bibitem{Kaj}
\newblock R. Kajikiya, T. Miyakawa,
\newblock On $L^2$ decay of weak solutions of the {N}avier-{S}tokes equations in $\mathbb{R}^n$,
\newblock  \emph{Math. Zeit.} \textbf{192} (1986), 135-148.

 \bibitem{Gr}
\newblock G. Karch, D. Pilarczyk,
\newblock Asymptotic stability of {L}andau solutions to {N}avier-{S}tokes system,
\newblock  \emph{Arch. Ration. Mech. Anal.} \textbf{202} (2011), 115-131.



\bibitem{Ka}
\newblock G. Karch, D. Pilarczyk, M.E. Schonbek,
\newblock  $L^2$-asymptotic stability of singular solutions to the Navier-Stokes system of equations in $\mathbb{R}^3$,
\newblock  \emph{J. Math. Pures Appl.} \textbf{108} (1) (2017), 14-40.



\bibitem{Kat} 
\newblock T. Kato,
\newblock Strong $L^p$-solutions of the Navier-Stokes equation in $\mathbb{R}^m$, with applications to weak solutions,
\newblock  \emph{Math. Zeit.} \textbf{187} (1984), 471-480.


\bibitem{Kik}
\newblock N. Kikuchi, G. Seregin,
\newblock Weak solutions to the Cauchy problem for the Navier-Stokes equations satisfying the local energy inequality,in Nonlinear equations and spectral theory
\newblock \emph{Amer. Math. Soc. Transl. Ser.}, vol. 2,  \textbf{220} (2007), 141-164.

 \bibitem{Ko} 
\newblock G. Koch, N. Nadirashvili, G. Seregin, V. \v{S}ver\'{a}k,
\newblock Liouville theorem for the {N}avier-{S}tokes equations and applications,
\newblock  \emph{ Acta Math.} \textbf{203} (1) (2009), 83-105.






 \bibitem{Kw} 
\newblock H. Kwon, T.-P.  Tsai,
\newblock Global {N}avier-{S}tokes flows for non-decaying initial
              data with slowly decaying oscillation,
\newblock  \emph{Comm. Math. Phys.} \textbf{375} (3) (2020), 1665-1715.


\bibitem{La}
\newblock L.D. Landau,
\newblock A new exact solution of the {N}avier-{S}tokes equations,
 \newblock  \emph{C. R. (Dokl.) Acad. Sci. URSS} \textbf{43} (1944), 286-288.


\bibitem{Lem1}
 \newblock P.G. Lemari\'{e}-Rieusset,
\newblock \emph{Recent Developments in the {N}avier-{S}tokes Problem},
\newblock Chapman \& Hall/CRC Research Notes in Mathematics, vol. 431, 2002.




\bibitem{Lem}
 \newblock P.G. Lemari\'{e}-Rieusset,
\newblock \emph{The Navier-Stokes Problem in the 21st Century},
\newblock CRC Press, Boca Raton, FL, 2016.




\bibitem{Ler} 
\newblock J. Leray,
\newblock Sur le mouvement d'un liquide visqueux emplissant l'espace,
\newblock \emph{ Acta Math.} \textbf{63} (1934), 193-248.


 \bibitem{Li1}
\newblock L. Li, Y.Y. Li, X. Yan,
\newblock Homogeneous solutions of stationary {N}avier-{S}tokes
             equations with isolated singularities on the unit sphere. {I}.
              {O}ne singularity,
 \newblock  \emph{Arch. Ration. Mech. Anal.} \textbf{227} (2018), 1091-1163.




 \bibitem{Li2}
\newblock L. Li, Y.Y. Li, X. Yan,
\newblock Homogeneous solutions of stationary Navier-{S}tokes
              equations with isolated singularities on the unit sphere.
              {II}. Classification of axisymmetric no-swirl solutions,
 \newblock  \emph{J. Differential Equations} \textbf{264} (2018), 6082-6108.




\bibitem{Li4}
\newblock L. Li, Y.Y. Li, X. Yan,
\newblock Vanishing viscosity limit for homogeneous axisymmetric
              no-swirl solutions of stationary {N}avier-{S}tokes equations,
 \newblock  \emph{J. Funct. Anal.} \textbf{277} (2019), 3599-3652.



 \bibitem{Li3}
\newblock L. Li, Y.Y.  Li, X. Yan,
\newblock Homogeneous solutions of stationary {N}avier-{S}tokes
              equations with isolated singularities on the unit sphere.
              {III}. Two singularities,
 \newblock  \emph{Discrete Contin. Dyn. Syst.} \textbf{39} (2019), 7163-7211.

\bibitem{Yan Li}
\newblock Y.Y.  Li, X. Yan,
\newblock Asymptotic stability of homogeneous solutions of incompressible stationary Navier-Stokes equations,
 \newblock arXiv:1911.03002.


\bibitem{Ma} 
\newblock K. Masuda,
\newblock Weak solutions of {N}avier-{S}tokes equations
 \newblock  \emph{Tohoku Math. J.} \textbf{36} (1984), 623-646.



  \bibitem{Ngu}
\newblock H.M  Nguyen, M. Squassina,
\newblock Logarithmic Sobolev inequality revisited,
 \newblock  \emph{C. R. Math.} \textbf{355}(4), (2017), 447-451.


\bibitem{Oga}
\newblock T. Ogawa, S. Rajopadhye, M. Schonbek,
\newblock Energy decay for a weak solution of the Navier-Stokes equation with slowly varying external forces,
 \newblock  \emph{J. Funct. Anal.} \textbf{144} (1997), 325-358.








\bibitem{Paz}
\newblock A. Pazy,
\newblock \emph{Semigroups of linear operators and applications to partial
              differential equations},
\newblock Springer-Verlag, New York, vol. 44, 1983.



\bibitem{Rob}
 \newblock J.C. Robinson, J.L. Rodrigo, W. Sadowski,
\newblock \emph{The three-dimensional {N}avier-{S}tokes equations},
\newblock Cambridge Studies in Advanced Mathematics, vol. 157, Cambridge University Press, Cambridge, Classical theory, 2016.

  \bibitem{Sch1} 
\newblock M.E. Schonbek,
\newblock Decay of solutions to parabolic conservation laws,
 \newblock  \emph{Commun. Partial Differ. Equ.} \textbf{7} (1980), 449-473.

  \bibitem{Sch} 
\newblock M.E. Schonbek,
\newblock $L^2$ decay for weak solutions of the {N}avier-{S}tokes
              equations,
 \newblock  \emph{Arch. Ration. Mech. Anal.} \textbf{88} (1985), 209-222.



\bibitem{Se}
\newblock G. Seregin, V. \v{S}ver\'{a}k,
\newblock On global weak solutions to the {C}auchy problem for the
              {N}avier-{S}tokes equations with large {$L_3$}-initial data,
 \newblock  \emph{Nonlinear Anal.} \textbf{154} (2017), 269-296.


\bibitem{Sl}
\newblock N.A. Slezkin,
\newblock On an integrability case of full differential equations of the motion of a
viscous fluid, in: Uchen. Zapiski  Moskov. Gosud. Universiteta, vol. 2, Gosud. Tehniko-Teoret. Izdat., Moskva/Leningrad, 1934, pp. 89-90.






\bibitem{St}
\newblock E.M. Stein,
\newblock \emph{Harmonic analysis: real-variable methods, orthogonality, and oscillatory integrals},
\newblock  Princeton University Press, Princeton, NJ, vol. 43, 1993.





\bibitem{Sv}
\newblock V. \v{S}ver\'{a}k,
\newblock On Landau’s solutions of the {N}avier-{S}tokes equations, Problems in mathematical analysis, No. 61,
 \newblock  \emph{J. Math. Sci. (NY)}, \textbf{179}(1) (2011), 208-228.



\bibitem{Swa}
\newblock  C. Swanson,
\newblock The best Sobolev constant,
 \newblock  \emph{Appl Anal.}, \textbf{47}(4) (1992), 227–239.





\bibitem{Te}
\newblock R. Temam,
\newblock \emph{Navier-Stokes Equations: Theory and Numerical Analysis},
\newblock  Reprint of the 1984 edition. AMS/Chelsea Publishing, Providence, 2001.

\bibitem{Ti}
\newblock G. Tian, Z.P. Xin,
\newblock One-point singular solutions to the {N}avier-{S}tokes equations,
 \newblock  \emph{Topol. Methods Nonlinear Anal.} \textbf{11} (1998), 135-145.




\bibitem{Tsa}
\newblock T.-P.  Tsai,
\newblock \emph{Lectures on {N}avier-{S}tokes equations},
\newblock American Mathematical Society, Providence, RI, vol. 192, 2018.








\bibitem{Wie} 
\newblock M. Wiegner,
\newblock  Decay results for weak solutions of the Navier-Stokes equations in $\mathbb{R}^n$,
\newblock  \emph{J. Lond. Math. Soc.} \textbf{35} (1987), 303-313.

\bibitem{Zhang}
\newblock J.J. Zhang, T. Zhang,
\newblock Local well-posedness of perturbed Navier-Stokes system around Landau solutions,
\newblock  \emph{preprint.}
\end{thebibliography}


\end{document}